\documentclass{amsart}
\usepackage{enumerate}
\usepackage{verbatim}
\usepackage{hyperref}
\usepackage{amsmath,amssymb,graphicx,float,epsf}
\usepackage{graphicx}
\usepackage{amsmath}
\usepackage{amsthm}
\usepackage{amsfonts}

\makeatletter
\newcommand{\labitem}[2]{%
\def\@itemlabel{\textbf{#1}}
\item
\def\@currentlabel{#1}\label{#2}}
\makeatother

\newcommand{\supp}{\mathrm{supp}}
\DeclareMathOperator*{\esssup}{ess\,sup}
\DeclareMathOperator*{\essinf}{ess\,inf}

\def\RR{\mathbb{R}}

\newtheorem{Def}{Definition}
\newtheorem{Lemma}{Lemma}
\newtheorem{Th}{Theorem}
\newtheorem{Cor}{Corollary}
\newtheorem{Prop}{Proposition}
\newtheorem{Rem}{Remark}
\allowdisplaybreaks
\numberwithin{equation}{section}

\begin{document}

\author[Korobenko]{Lyudmila Korobenko}
\address{University of Calgary\\
Calgary, Alberta\\
lkoroben@ucalgary.ca}
\author[Rios]{Cristian Rios}
\address{University of Calgary\\
Calgary, Alberta\\
crios@ucalgary.ca}
\title[Regularity of Non-Doubling Second Order Equations]{Regularity of Solutions to Degenerate Non-Doubling Second Order Equations}
\thanks{The second author is supported by the Natural Sciences and
Engineering Research Council of Canada.}

\begin{abstract}
We prove that every weak solution to a certain class of infinitely degenerate quasilinear equations is continuous. An essential feature of the operators we consider is that their Fefferman-Phong associated metric may be non doubling with respect to Lebesgue measure. 
\end{abstract}

\maketitle

\section{Introduction}

One of the fundamental results of the theory of second order elliptic equations is the De Giorgi-Nash-Moser a-priori regularity of weak solutions \cite{DeGiorgi:1957, Nash:1958, Moser:1960, Moser:1961}, see also \cite{Gilb}.  Given a second order operator $L$ of the form 
\begin{equation}\label{L_linear}
L=\nabla^TA(x)\nabla
\end{equation}
satisfying an ellipticity condition
\begin{equation}\label{ellipticity}
\lambda\, |\xi|^2 \leq \xi^TA(x)\xi\leq \Lambda\,|\xi|^2,\qquad\forall \xi\in\mathbb{R}^n\qquad\hbox{and some } \Lambda>\lambda>0,
\end{equation}
a weak solution to the equation $Lu=f$ in an open bounded domain $\Omega\subset\RR^n$, with $f\in L^{\infty}(\Omega)$, is any function $u$ from the Sobolev space $W^{1,2}(\Omega)$ such that
$$
-\int(\nabla u)^TA\nabla w=\int fw
$$
for all $w\in W^{1,2}_0(\Omega)$. The classical theory tells us that weak solutions must be H\"{o}lder continuous, i.e. $u\in C^{\alpha}(\Omega)$ for some $0<\alpha\le 1$ which only depends on $n$, the ellipticity constant $\lambda$, and the upper bound $\Lambda$. In case of a quasilinear operator, i.e. when the matrix of coefficients $A$ itself depends on a solution, $A=A(x,u)$, the ellipticity condition can be defined in a similar way, and under some structural assumptions the Moser techniques can be extended to this type of operators to obtain H\"{o}lder continuity of solutions \cite{Lady68, Gilb}. The techniques pioneered by De Giorgi, Nash, and Moser have an extensive and prolific evolution: \cite{Morrey61, Serrin63, KrSa80, Tr83, FKS82, FG85, Saw, Ma13}, among many other related results. In virtually all of this existing theory on a-priori regularity of weak solutions (to equations of elliptic or degenerate elliptic type), the space geometry is of ``homogeneous type''. This means that the topology is given by a metric or quasi-metric for which the underlying measure is doubling, that is, the measure of a ball of radius $2r$ is bounded by a fixed multiple of the measure of the ball of radius $r$ with the same center. This is the case for all works just cited. 

A wider class than elliptic operators is that of subelliptic operators, which are defined by differential inequalities of the type
\[\| u\|_\varepsilon^2\le C\left( \langle Lu,u\rangle+\| u\|^2\right), \]
where $\varepsilon>0$, $\|\cdot\|$ is the $L^2$-norm, $\|\cdot\|_\varepsilon$ is the Sobolev norm, and $\langle\cdot,\cdot\rangle$ is the inner product in $L^2$. These inequalities say that weak solutions have some a-priori ``gain'' in regularity.  When the operators are given in terms of sum of squares of vector fields, subellipticity has been characterized by the celebrated H\"{o}rmander condition \cite{Hor}, requiring the Lie group generated by the vector fields to be of finite type. In the case of nonnegative definite linear of second order, Fefferman and Phong \cite{Fef} obtained an enlightening characterization of subellipticity in terms of their subunit metric balls $B_L$ and the Euclidean balls $B_E$: an operator $L$ is subelliptic if and only if there exist positive constants $C$ and $\epsilon$ such that 
\begin{equation}\label{FP} 
B_E(x,r)\subseteq B_L(x,Cr^\epsilon)
\end{equation}
for all $x$ and $r>0$ in the domain of consideration. What is remarkable about Fefferman and Phong's result is that it makes an apparent connection between the geometry inherent to the operator and the regularity of weak solutions. Note that their condition \eqref{FP} implies that the metric balls $B_L$ are doubling with respect to Lebesgue's measure. 

In \cite{Saw} Sawyer and Wheeden considered quasilinear equations in divergence form, and defined notions of subelliptic operators in terms of a priori H\"{o}lder regularity of weak solutions among classes of operators with either bounded coefficients or, more generally, with $L^p$ coefficients. Some of the main results in \cite{Saw} are a generalization of the H\"{o}rmander criterium for operators with rough coefficients, and a generalization of the Fefferman-Phong condition for operators with rough coefficients. The proper definition of weak solutions to quasilinear equations requires Sobolev spaces adapted to the nonnegative definite form induced by, or controlling, the principal coefficients of the operator. These authors also developed in detail the generalized theory of Sobolev spaces necessary to treat these more general operators \cite{Saw2}, we use their definitions in our present work. The main elements in \cite{Saw} to guarantee subellipticity for operators with rough coefficients are:
\begin{enumerate}[(i)]
\item \label{itm-doub} the doubling condition of the metric balls, 
\item \label{itm-cont} a containment condition of the type \eqref{FP},
\item \label{itm-sobo} a Sobolev inequality
\[ \left\{ \frac{1}{|B|}\int_B |w|^{2\sigma}\right\}^{\frac{1}{2\sigma}} \le C r \left\{ \frac{1}{|B|}\int_B [\nabla w]^{2}_{Q}\right\}^{\frac{1}{2}} +C \left\{ \frac{1}{|B|}\int_B |w|^{2}\right\}^{\frac{1}{2}} ,\]
\item \label{itm-poin} a Poincar\'e inequality
\[ \left\{ \frac{1}{|B|}\int_B \left| w-\left( \frac{1}{|B|}\int_Q w\right)\right|^{2\sigma}\right\}^{\frac{1}{2\sigma}} \le C r \left\{ \frac{1}{|B^*|}\int_{B^*} [\nabla w]^{2}_{Q}\right\}^{\frac{1}{2}} ,\]
\item \label{itm-aslc} ``accumulating sequence of Lipschitz cutoff functions'': there exist positive constants $C_0$ and $N$ such that for each ball $B$ of radius $r$ there is a sequence of Lipschitz functions $\{\psi_j\}_{j=1}^\infty$ on $B$ with the following properties
\[ \left\{ \begin{array}{ll} 
         \textrm{supp}(\psi_j)& \subset B,\\ 
        B^* & \subset \{ x:\psi_j(x)=1\},\ j\ge 1\\ 
        \textrm{supp}(\psi_{j+1}) &  \Subset\{ x:\psi_j(x)=1\},\ j\ge 1\\
        \left\{ \frac{1}{|B|}\int [\nabla\psi_j]_Q^pdx\right\}^{\frac{1}{p}} & \le C_p\frac{j^N}{r}, j\ge 1,
    \end{array} \right.\]
\end{enumerate}
Where $B$ is a ball of radius $r$, $B^*$ is a ball with the same center and radius $C_0r$, and $[\nabla w]_Q^2$ is given by \eqref{eq-qgrad}. 

One of the first systematic approaches to a-priori regularity for more general elliptic operators is found in \cite{FKS82}, where linear operators in divergence form were considered, with ellipticity controlled by an $A_2$ Muckemphout weight. In that paper it is shown that if properties \eqref{itm-doub}--\eqref{itm-poin} are satisfied, plus a condition on uniqueness of gradients in the weighted Sobolev spaces, then positive solutions satisfy Harnack's inequality and consequently weak solutions are H\"{o}lder continuous. Since the geometry considered in \cite{FKS82} is Euclidean, property \eqref{itm-aslc} also trivially holds. It was later shown by other authors that property \eqref{itm-sobo} and the uniqueness condition on the gradient are consequences of properties \eqref{itm-doub} and \eqref{itm-poin} \cite{Saw2, SC, Haj2, HK95}. In a more recent development \cite{Kor}, the present authors and Maldonado proved that properties \eqref{itm-sobo} and \eqref{itm-aslc} imply property \eqref{itm-doub}, the doubling condition.

In the classical work of H\"{o}rmander \cite{Hor} it is shown that subellipticity implies hypoellipticity. An operator is hypoelliptic if the coefficients and the right-hand side of the equation are $C^{\infty}$, then the solution must also be $C^{\infty}$. The class of hypoelliptic operators is wider than the class of subelliptic operators \cite{Fedii, Kohn}. It is enlightening to study the relation between the order to which ellipticity fails, i.e. the order at which eigenvalues $\lambda(x)$ of the coefficients matrix become zero, and the a-priori regularity of solutions. First, if $\lambda(x)\equiv 0$ on a set of positive measure, then the operator can fail to be hypoelliptic \cite[Example 48]{Saw}. 
On the other hand, if $\lambda(x)$ vanishes as a polynomial (finite degeneracy) then the operator is subelliptic \cite{Rios}. The intermediate case of infinite degeneracy of eigenvalues, for example, when $\lambda(x)$ vanishes together with all its derivatives on hyperplanes, has been considered by many authors \cite{Kus,Mor,Christ}. A number of important regularity results have been established for different classes of linear infinitely degenerate elliptic operators \cite{Fedii, Kohn, Mor, Tri}. An analysis of quasilinear degenerate operators is much more complicated and the theory for this type of operators is less developed.

In our present work we build on the approach in \cite{Fef, Saw, Saw2} and consider operators with coefficients with very little regularity, what makes it possible to include applications to quasilinear operators. More importantly, we allow the underlying metric to be non-doubling. The type of degeneracy assumed is quite general, and it allows for infinite vanishing of the coefficients, what excludes subellipticity.  We consider a second order quasilinear equation of the form 
\begin{equation}\label{eq_0}
Lu:=\nabla'A(x,u(x))\nabla u=f
\end{equation} 
and assume structural conditions
\begin{equation}\label{struc_0}
k\, \xi^TQ(x)\xi \leq \xi^TA(x,z)\xi\leq K\,\xi^TQ(x)\xi
\end{equation}
on the quasilinear matrix $A(x,u(x))$. To state our main result about the regularity of weak solutions we assume the existence of a certain metric $d$ for which the metric balls $B(x,r)=\{ y\in\Omega:\ d(x,y)<r\}$ define mutually equivalent topologies with the Euclidean metric. That is, for all $x\in\Omega$ and $r>0$ the metric balls $B(x,r)$ are open sets in the Euclidean sense, and the Euclidean balls $B_E(x,r)=\{y:|x-y|<r\}$ which are contained in $\Omega$ are open sets with respect to the metric $d$.  For simplicity we will further assume that the metric balls are uniformly bounded in the Euclidean distance, in the sense stated in the containment condition \eqref{upper} below. We state our assumptions more precisely in the form of the following two conditions:
\begin{itemize}
\item $\forall x\in\Omega\;\exists\alpha_x:\RR_{+}\rightarrow\RR_{+}$ increasing and such that $\alpha_x(r)>0$ for any $r>0$ and
\begin{equation}\label{contain}
B_E(x,\alpha_x(r))\subset B(x,r),\  \forall x\in\Omega,\  \forall r>0\;\hbox{s.t.}\;B(x,r)\subset\Omega
\end{equation}
\item There exists a constant $C>0$ s.t. for any $x\in\Omega$, any $r>0$ s.t. $E(x,Cr)\subset\Omega$ there holds
\begin{equation}\label{upper}
B(x,r)\subset B_E(x,Cr).
\end{equation}
\end{itemize}
Equivalently, for every $x\in\Omega$ the functions $d_x(y)=d(x,y)$ and $|y|_x=|x-y|$ are continuous in the Euclidean topology and in the $d$-topology, respectively. Condition \eqref{upper} says that the modulus of continuity of $|\cdot|_x$ is of the form $Cd_x(\cdot)$; in particular, for any $x\in\Omega$ there exists $\eta>0$ such that $B(x,r)\subset\Omega$ for any $0<r<\eta\, dist(x,\partial\Omega)$.

A natural choice of the metric $d$ is the subunit metric which will be discussed in Section \ref{Subunit metric}. All of the results, however, are axiomatic in the sense that if one can find any metric satisfying the required assumptions the results will hold true. For other choices of metric balls, e.g. flag balls and adapted noninterference balls, see \cite{Saw}.\\

We will require the following Sobolev inequality:
there exist $\sigma>1$ and $\eta>0$ such that for all balls $B=B(y,r)$ with $y\in\Omega,\  0<r<\eta\, dist(y,\partial\Omega)$ there holds
\begin{equation}\label{sobolev}
\begin{split}
\left\{\frac{1}{|\supp(w)|}\int\limits_{B}|w|^{2\sigma}\right\}^{\frac{1}{2\sigma}}\leq &\ Cr \,\left\{\frac{1}{|\supp(w)|}\int\limits_{B}[\nabla w]_{Q}^2\right\}^{\frac{1}{2}}\\
&+C\left\{\frac{1}{|\supp(w)|}\int\limits_{B}|w|^{2}\right\}^{\frac{1}{2}}
\end{split}
\end{equation}
for all $(w,\nabla w)\in\left(\mathcal{W}_{Q}^{1,2}\right)_0(B)$, the closure in $\mathcal{W}_{Q}^{1,2}(B)$ of $(w,\nabla w)$ where $w\in Lip_c(B)$. The spaces $\mathcal{W}_{Q}^{1,p}$ are the strong degenerate Sobolev spaces associated to the matrix $Q(x)$ as defined in Section \ref{weak solutions}. The gradient on the right-hand side is the $Q$-gradient, defined in the standard way
\begin{equation}\label{eq-qgrad}
[\nabla \varphi]_{Q}^2:=\nabla \varphi^TQ\nabla \varphi.
\end{equation}
Note that it is slightly different from a classical Sobolev inequality \eqref{itm-sobo} used, for example, in \cite{Saw,Wheeden}. Here we are taking ``true averages'' of the function $w$ by dividing by the measure of the support, and not the measure of any ball containing the support. The reason is that in the non doubling case by varying the size of the ball it might be possible to strengthen the inequality. Moreover, as we show in \cite{Kor}, the classical Sobolev inequality together with the accumulating sequence of Lipschitz cutoff function introduced in \cite{Saw} in fact imply the doubling condition on the metric balls. Therefore, in order to be able to include a non doubling case, the classical Sobolev inequality must be weakened, and we believe (\ref{sobolev}) is the right version. We also note that in the doubling case the two versions coincide.

We also require the following strong Poincar\'{e} inequality: there exists $\eta>0$ such that for all balls $B=B(y,r)$ with $y\in\Omega,\  0<r<\eta\, dist(y,\partial\Omega)$ there holds
\begin{equation}\label{poincare}
\int\limits_{B}\left|w-\left\langle w\right\rangle_B\right|\leq Cr\int\limits_{ B}[\nabla w]_{Q}
\end{equation}
for all $(w,\nabla w)\in\mathcal{W}_{Q}^{1,2}(B)$ where we denote $\left\langle w\right\rangle_B=\left(1/|B|\right)\int_{B}w$.

Finally, we also assume a growth condition on the non doubling order of the metric balls, and the existence of an accumulating sequence of Lipschitz cutoff functions.

\begin{Th}\label{main_intro}
Suppose that $A(x,z)$ is a nonnegative semidefinite matrix in $\Omega\times\mathbb{R}$ and it satisfies (\ref{struc_0}). Let $d(x,y)$ be a symmetric metric in $\Omega$, and $B(x,r)=\{y\in\Omega:d(x,y)<r\}$ with $x\in\Omega$ are the corresponding metric balls. Then every weak solution of (\ref{eq_0}) is continuous provided that
\begin{enumerate}
\item the containment condition (\ref{contain}) holds, 
\item the boundedness condition (\ref{upper}) holds,
\item the Sobolev and Poincar\'{e} inequalities (\ref{sobolev}) and (\ref{poincare}) hold,
\item the non doubling order $\delta_x(r)$ of metric balls, defined by \eqref{nondoub_order}, satisfies the growth condition \eqref{nondoub_order_control} for every $x\in\Omega$, and
\item there exist accumulating sequences of cutoff functions satisfying \eqref{cutoff} and \eqref{spec_cutoff}.  
\end{enumerate}
\end{Th}
First, as we mentioned earlier, a big strength of this result is that it does not rely on the doubling assumption for metric balls. We cannot, however, allow for all types of non doubling balls and we must limit the blow up rate of the ratio $|B(x,2r)|/|B(x,r)|$ as $r\rightarrow 0$. On the other hand, since our containment condition is much weaker than the one used by Fefferman and Phong \cite{Fef} and by Sawyer and Wheeden \cite{Saw}, the best possible a-priori regularity is that weak solutions are continuous and not, in general, H\"{o}lder continuous. Our containment condition \eqref{contain} only requires equivalence of topologies and boundedness of the metric balls. The Fefferman-Phong containment condition requires that every subunit ball contains a Euclidean ball with radius proportional to a power of the original radius; this requirement can be loosely thought of as a condition of finite degeneracy, while our condition allows for infinite degeneracy. 

Secondly, the requirement of an accumulating sequence of cutoff functions is different from the one in \cite[Theorem 8]{Saw} but it is shown to hold for subunit metric spaces under mild restrictions on the operator, see Lemma \ref{spec_cutoff_lemma} in Section \ref{section_cutoff_functions}. The main difference is that in the non doubling setting, the measures of the supports of cutoff functions ``accumulate'' when performing the Moser iteration. The supports therefore must be chosen in a very special way to obtain convergence of the bounds. The Sobolev and Poincar\'{e} inequalities are much more difficult to establish in the setting of non doubling metric measure spaces, there are still many open problems related to them, and they are subjects for further research.

Finally, note that this result only gives the continuity of weak solution, and not the full hypoellipticity result. However, it has been shown by Rios et al. \cite{Rios2} that every weak continuous solution to a certain class of infinitely degenerate quasilinear equations is smooth provided the coefficients are smooth. Therefore, the theorem may provide a final important bridge between weak solutions and hypoellipticity for some quasilinear equations.

The paper is organized as follows. We first introduce the concept of a non doubling order and give related definitions in Section \ref{Non doubling metrics}. This is precisely the quantity that we need to control in order to obtain continuity of weak solutions. Next, in Section \ref{weak solutions} we introduce degenerate Sobolev spaces needed to define weak solutions to degenerate second order equations. This allows to give the widest possible definition of weak solutions.
Section \ref{Proof of the main result} is dedicated to the proof of our main result, the continuity of weak solutions. It is then shown in Section \ref{Subunit metric} that some of the requirements of Theorem \ref{main_intro} are satisfied by subunit metric balls under certain conditions on the operator. Finally, Section \ref{examples} explores some examples of linear and quasilinear operators that fall under the framework of Theorem \ref{main_intro}. In these examples an operator has an infinite degeneracy at the origin, and the weak containment condition is shown to hold for subunit metric balls associated to this operator.

\section{Non doubling metrics}\label{Non doubling metrics}

We now give definitions and discuss some properties of metric balls that are non doubling with respect to Lebesgue measure. Let $\Omega\subset\RR^n$ be an open bounded domain and $d$ a metric on $\Omega$.

\begin{Def}
Metric balls are said to satisfy the doubling condition with respect to the measure $\mu$ if
$$
\mu(B(x,2r))\leq C\mu(B(x,r)),\quad x\in\Omega,\,0<r<\infty
$$
In this case a metric measure space $(\Omega,d,\mu)$ is called a space of homogeneous type.
\end{Def}

Spaces of homogeneous type are of particular interest in applications to differential equations because a great portion of useful results from Euclidean spaces can be extended to space of homogeneous type. These include Sobolev spaces on metric spaces \cite{Shan, Haj}, BMO spaces \cite{Tolsa}, and singular integral operators \cite{Naz-Treil-Volb}. 

Even in the doubling case, in order to establish regularity results, one might have to impose certain conditions on the speed of growth/decay of the metric balls. Examples of such conditions are $r^s\leq C|B(x,r)|$ and $|B(x,r)|\leq Cr^N$. Our approach is to require some control on the ratio $|B(x,2r)|/|B(x,r)|$ as $r\to 0$, yet we still allow for rates smaller than any power so the metric may be non-doubling. 

\begin{Def}
Let $\Omega\subset\RR^n$ and $d$ a metric on $\Omega$. For $x\in\Omega$ we say that the metric balls are non doubling of order $\delta_x(r)$ at $x$ if there exist constants $C>1$ and $\eta>0$ such that
\begin{equation}\label{nondoub_order}
5/4\leq\frac{|B(x,r+\delta_x(r))|}{|B(x,r)|}\leq 2C,\quad \forall 0<r<\eta\, dist(x,\partial\Omega)
\end{equation}
and $\delta_x(r)/r\to 0$ as $r\to 0$. Note, that $\delta_x(r)$ is not uniquely determined by (\ref{nondoub_order}) but by the order of non doubling we will implicitly understand the smallest such $\delta$, which must be positive for $r>0$.
\end{Def}
We will sometimes omit the dependence on $x$, when the reference point is clear, and write $\delta(r)$ for $\delta_x(r)$. Note, that if $\delta_x(r)\geq Cr$ for all $r<\eta\, dist(x,\partial\Omega)$ then the metric balls are in fact doubling. Otherwise, we have $\lim_{r\to 0}(\delta_x(r)/r)=0$. It is also reasonable to assume that $\delta_x(r)$ is an increasing function of $r$. The rate of vanishing of the ratio $\delta_x(r)/r$ plays an important role, for example, an exponential rate, $\delta_x(r)/r\approx \exp(-1/r^2)$ is much harder to handle than the linear rate $\delta_x(r)/r\approx r$. It is easy to see that the non doubling order $\delta_x(r)$ can be connected to the ``doubling ratio'' $|B(x,2r)|/|B(x,r)|$. Indeed, let $m$ be the first integer such that $2r-m\delta_x(r)\leq r$ and $2r-(m-1)\delta_x(r)>r$, then
\begin{equation*}
\begin{split}
\frac{|B(x,2r)|}{|B(x,r)|}&= \frac{|B(x,(2r-\delta_x(r))+\delta_x(r))|}{|B(x,2r-\delta_x(r))|}\cdot\frac{|B(x,(2r-2\delta_x(r))+\delta_x(r))|}{|B(x,2r-2\delta_x(r))|}\cdot\ldots\\
&\ldots\cdot\frac{|B(x,(2r-m\delta_x(r))+\delta_x(r))|}{|B(x,r)|}\leq (2C)^m\leq (2C)^{\frac{r}{\delta_x(r)}+1}
\end{split}
\end{equation*}
with $C$ as in (\ref{nondoub_order}). Similarly, it can be shown that
$$
\frac{|B(x,2r)|}{|B(x,r)|}\geq \left(\frac{5}{4}\right)^{\frac{r}{\delta_x(2r)}-1}.
$$
Therefore, we might require the bounds on either the doubling ratio $|B(x,2r)|/|B(x,r)|$ or the non doubling order $\delta_x(r)$. For our analysis we have chosen to do the latter, see the proof of Theorem \ref{main_intro} in Section \ref{Proof of the main result}.

The next condition on $\delta_x(r)$ will guarantee the continuity of weak solutions (see Theorem \ref{continuity_thm}) for certain values of the parameters $\lambda$ and $C$
\begin{equation}\label{nondoub_order_control}
\ln r \ln\left(1-\frac{\exp(-(r/\delta(r))^{\lambda})}{2C}\right)\to \infty,\quad \hbox{as}\quad r\to 0
\end{equation}
for all $x\in\Omega$, $r>0$ small enough, $\lambda>1$ and $C>1$. Note, that the above condition implies a very slow vanishing of $\delta_x(r)/r$, in particular, it must be slower than any power of $r$. Some model examples are given in Section \ref{examples}.

\subsection{Cutoff Functions}
We can now formulate the requirements of the existence of certain cutoff functions employed later in the proofs. The accumulating sequence is an adaptation of the one introduced by Sawyer and Wheeden \cite{Saw}, \cite{Saw2}, to the case of non doubling metric balls. We assume there are positive constants $\nu$, $N$ and $\eta$ such that for each ball $B(y,r)$ with $y\in\Omega$, $0<r<\eta\, dist(y,\partial\Omega)$, there is a sequence of Lipschitz cutoff functions $\left\{\psi_j\right\}_{j=1}^{\infty}$ with the following properties:
\begin{equation}\label{cutoff}
\begin{cases}
E_1&=\supp(\psi_1)\subset B(y,r),\\
B(y,\nu r) &\subset \{x:\psi_j(x)=1\},\;j\geq 1,\\
E_{j+1}=\supp(\psi_{j+1}) &\Subset \{x:\psi_j(x)=1\},\;j\geq 1,\\
\frac{\displaystyle |E_j|}{\displaystyle|E_{j+1}|}&\leq C,\;j\geq 1,\\
\psi_j &\hbox{is Lipschitz},\;j\geq 1,\\
||[\nabla\psi_j]_Q||_{L^{\infty}(B(y,r))} &\leq \frac{C}{\displaystyle (1-\nu)\delta(\nu r)\left(1-\frac{\delta(\nu r)}{r}\right)^{j}},\;j\geq 1,
\end{cases}
\end{equation}
where we denote $[\nabla\psi_j]_Q=((\nabla\psi_j)^T Q\nabla\psi_j)^{1/2}$. This sequence is used in the process of implementing the Moser iteration in Propositon \ref{Moser iteration} Section \ref{Proof of the main result}. For special cases of subunit metrics associated to the operator $L$, the existence of such sequence is shown in Section \ref{section_cutoff_functions}. We also require the existence of a Lipschitz cutoff function satisfying the following condition
\begin{equation}\label{spec_cutoff}
\begin{cases}
\supp(\phi_{r}) &\subseteq B(y,r+\delta(r)),\\
\{x:\phi_{r}(x)=1\} &\supseteq B(y,r+\delta(r)/2),\\
||[\nabla\phi_{r}]_Q||_{L^{\infty}(B(y,r)} &\leq \frac{C}{\delta(r)},
\end{cases}
\end{equation}
where again $[\nabla\phi_{r}]_Q=((\nabla\phi_{r})^T Q\nabla\phi_{r})^{1/2}$. This cutoff function is used in the proof of weak logarithmic estimates, Lemma \ref{lemma_weak1}. 

\section {Degenerate Sobolev spaces and weak solutions}\label{weak solutions}

Let $\Omega\subset\mathbb{R}^n$ be an open bounded domain. Consider the following second order quasilinear equation
\begin{equation}\label{eq}
Lu=\nabla^TA(x,u(x))\nabla u=f
\end{equation}
in $\Omega$ where $f\in L^{\infty}(\Omega)$ and the matrix $A=\{a_{ij}\}$ is nonnegative semidefinite and can degenerate to infinite order. We will assume that there exists a nonnegative semidefinite locally integrable matrix $Q(x)$ and constants $K\geq k>0$ such that $A(x,z)$ satisfies the following structural condition for a.e. $x\in\Omega$ and all $z\in\mathbb{R}$, $\xi\in\mathbb{R}^n$:
\begin{equation}\label{struc}
k\, \xi^TQ(x)\xi \leq \xi^TA(x,z)\xi\leq K\,\xi^TQ(x)\xi .
\end{equation}
Note that if we denote $\tilde{A}(x)=A(x,u(x))$ for a particular solution $u(x)$ then (\ref{struc}) obviously holds with $\tilde{A}(x)$ in place of $A(x,z)$.
To define solutions to (\ref{eq}) in the weakest sense possible we adopt the notion of strong degenerate Sobolev spaces \cite{Saw2}. 

\begin{Def}\cite{Saw2}
A form-weighted vector-valued $L^2$-space $\mathcal{L}^2(\Omega,Q)$ is a space consisting of all measurable $\RR^n$-valued functions $\textbf{f}(x),\;x\in\Omega$, satisfying
\begin{equation*}
||\textbf{f}||_{\mathcal{L}^2(\Omega,Q)}=\left\{\int_{\Omega}\textbf{f}(x)^TQ(x)\textbf{f}(x)dx\right\}^{\frac{1}{2}}
=\left\{\int_{\Omega}[\textbf{f}(x)]^2_Q dx\right\}^{\frac{1}{2}}<\infty
\end{equation*}
where we denote $[\textbf{U}(x)]^{2}_{Q}=\textbf{U}(x)^TQ(x)\textbf{U}(x)$ for any vector-valued function $\textbf{U}(x)$.
\end{Def}

As usual, we identify measurable $\mathbb{R}^n$-valued functions $\textbf{f}$ and $\textbf{g}$ satisfying $||\textbf{f}-\textbf{g}||_{\mathcal{L}^2(\Omega,Q)}=0$. We then denote by $\mathcal{L}^2(\Omega,Q)$ the space of equivalence classes of measurable $\mathbb{R}^n$-valued functions. It has been shown \cite{Saw2} that $\mathcal{L}^2(\Omega,Q)$ is complete with respect to the associated norm, and moreover is a Hilbert space with respect to the inner product
\begin{equation*}
\left\langle\textbf{f},\textbf{g} \right\rangle_{\mathcal{L}^2(\Omega,Q)}=\int_{\Omega}\textbf{f}(x)^TQ(x)\textbf{g}(x)dx.
\end{equation*}

\begin{Def}\cite{Saw2}
The space $W_{Q}^{1,2}(\Omega)$ is defined to be the completion of $Lip(\Omega)$ under the norm
\begin{equation*}
||w||_{W_{Q}^{1,2}(\Omega)}=\left\{\int_{\Omega}(|w|^2+[\nabla w]^{2}_{Q})\right\}^{\frac{1}{2}}
\approx ||w||_{L^2(\Omega)}+||\nabla w||_{\mathcal{L}^2(\Omega, Q)}.
\end{equation*}
By $\left(W_{Q}^{1,2}\right)_0(\Omega)$ we mean the closure of $Lip_c(\Omega)$ in $W_{Q}^{1,2}(\Omega)$ where $Lip_c(\Omega)$ is the space of Lipschitz functions with common compact support.
\end{Def}

The space $W_{Q}^{1,2}(\Omega)$ is thus a Banach space of Cauchy sequences in $Lip(\Omega)$. If $\{w_k\}_{k=1}^{\infty}$ is a Cauchy sequence in $W_{Q}^{1,2}(\Omega)$, then there are elements $w\in L^2(\Omega)$ and $\textbf{v}\in\mathcal{L}^2(\Omega,Q)$ such that $w_k\to w$ in $w\in L^2(\Omega)$ and $\nabla w_k\to\textbf{v}$ in $\mathcal{L}^2(\Omega,Q)$. The pair $(w,\textbf{v})$ represents the equivalence class containing the Cauchy sequence $\{w_k\}_{k=1}^{\infty}\in W_{Q}^{1,2}(\Omega)$ and we write $(w,\textbf{v})\in\mathcal{W}^2(\Omega, Q)$.
It is clear, however, that if $Q(x)$ is degenerate, the function $\textbf{v}\in\mathcal{L}^2(\Omega, Q)$ is not uniquely determined by $w\in W_{Q}^{1,2}(\Omega)$, see, for example, \cite[Section 4.3]{Rodney}. Therefore, $\nabla w$ denotes one of such vector-valued functions in $\mathcal{L}^2(\Omega, Q)$. We write $(w,\nabla w)\in \mathcal{W}_{Q}^{1,2}(\Omega)$ for $w\in W_{Q}^{1,2}(\Omega)$. An element $(w,\nabla w)\in \mathcal{W}_{Q}^{1,2}(\Omega)$ is said to be nonnegative if $w\geq 0$.

We now define weak solutions to (\ref{eq}) using the degenerate Sobolev spaces introduced above.

\begin{Def}\label{sol_def}
A pair $(u,\nabla u)\in \mathcal{W}^{1,2}_{Q}(\Omega)$ is a weak solution of (\ref{eq}) in $\Omega$ if
\begin{equation}\label{sol}
-\int(\nabla u)^TA\nabla w=\int fw
\end{equation}
for every nonngegative $w\in \left(W_{Q}^{1,2}\right)_0(\Omega)$. A function $u\in W^{1,2}_{Q}(\Omega)$ is called a weak subsolution (supersolution) if the above holds with $\geq(\leq)$ in place of equality.
\end{Def}

\subsection{Sub- and super- solutions}

Using equation (\ref{sol}) it is possible to find equations satisfied by a nonlinear function of $u$. These equations will be used in the process of performing Moser iterations --- one of the main steps in the proof of continuity. Typical examples of nonlinear functions used in Moser iterations are power functions $u^\beta$. We need to be careful, however, since sometimes it turns out to be necessary to truncate these functions. Moreover, we will also need to consider certain logarithmic functions of weak (sub-, super-) solutions. Therefore, it will be convenient to consider a class of ``admissible'' nonlinear functions \cite{Saw} that we will further compose with $u\in W_{Q}^{1,2}(\Omega)$.

\begin{Def}
Let $I\in\mathbb{R}$ be an interval and $h\in C^1(I)\cap C^2_{pw}(I)$ be positive and monotone, where $C^2_{pw}(I)$ is the space of piecewise twice continuously differentiable functions on $I$. The function $h$ is said to be admissible on $I$ if there exists a positive constant $C$ such that
$$
|h'(t)|,\;|h''(t)|,\;|th''(t)|\leq C,\quad t\in I.
$$
Moreover, given $u\in W_{Q}^{1,2}(\Omega)$, we say that $h$ is admissible for $u$ if $h$ is admissible on some interval $I$ containing the range of $u$.
\end{Def}

When composing $h''\in L^{\infty}$ with $u\in W_{Q}^{1,2}(\Omega)$ we might run into problems when $u$ takes values in the set of discontinuities of $h''$. To make sense of the expressions like $h''(u)\nabla u$ we will use the following proposition \cite{Saw2}
%

\begin{Prop}\cite[Proposition 22]{Saw2}
Suppose that $(u,\nabla u)\in \mathcal{W}_{Q}^{1,2}(\Omega)$ where $\Omega$ is bounded, and let
$$
\mathcal{R}_u=\{\alpha\in\mathbb{R}:u=\alpha\;\hbox{on a set of positive measure}\}.
$$
With $\nabla_{reg}u=\chi_{\{x\in\Omega:u(x)\notin\mathcal{R}_u\}}\nabla u$ we have $(u,\nabla_{reg} u)\in \mathcal{W}_{Q}^{1,2}(\Omega)$ and $(u,\nabla_{reg} u)$ satisfies
\begin{equation}\label{regular}
\|\chi_{\{u=\alpha\}}\nabla_{reg} u\|_{\mathcal{L}^2(\Omega, Q)}=0\quad\forall\alpha\in\mathbb{R}.
\end{equation}
\end{Prop}

We assume that all the elements of e adopt the following convention \cite{Saw2}
\begin{enumerate}[(I)]
\item  \label{itm-reggrad} If $u\in W_{Q}^{1,2}(\Omega)$ then $h(u)$ refers to the pair $(h(u),h'(u)\nabla_{reg}u)\in \mathcal{W}_{Q}^{1,2}(\Omega)$,
\item \label{itm-product-rule} if $u,v\in W_{Q}^{1,2}(\Omega)$ are represented by the pairs $(u,\nabla u),(v,\nabla v)\in \mathcal{W}_{Q}^{1,2}(\Omega)$ then the product $uv$ is represented by the pair $(uv,u\nabla v+v\nabla u)\in \mathcal{W}_{Q}^{1,1}(\Omega)$.
\end{enumerate}

\begin{Lemma}\cite[Lemma 19]{Saw2}\label{conv_prop}
Let $(u,\nabla u)\in \mathcal{W}_{Q}^{1,2}(\Omega)$ where $\Omega$ is bounded. If $f\in C^1(\mathbb{R})$ with $f'\in L^{\infty}(\mathbb{R})$ then $(f(u),f'(u)\nabla u)\in\mathcal{W}_{Q}^{1,2}(\Omega)$.
\end{Lemma}
\begin{Cor}\label{h_bound}
If $u\in W_{Q}^{1,2}(\Omega)$ and $h$ is admissible for $u$, it follows that $\tilde{u}=h(u)$ belongs to $W_{Q}^{1,2}(\Omega)$ provided $\Omega$ is bounded.
\end{Cor}

Any admissible nonlinear function of a weak (sub-, super-) solution to (\ref{eq}) satisfies a related equation in the weak sense, but for a smaller class of test functions. Following \cite{Saw}, \cite{Saw2} we now introduce weaker classes of solutions to (\ref{eq}).

\begin{Def}\label{w_solutions}
Let $\mathcal{W}$ be a subset of nonnegative elements in $\left(W_{Q}^{1,2}\right)_0(\Omega)$. A function $u\in W_{Q}^{1,2}(\Omega)$ is a $\mathcal{W}$-weak solution (subsolution, supersolution) to (\ref{eq}) in $\Omega$ if the integrals in (\ref{sol}) are absolutely convergent and the indicated (in)equality holds for all $w\in\mathcal{W}$.
\end{Def}
We can now define the classes of ``admissible test functions'' which will be used as $\mathcal{W}$-classes in the Definition \ref{w_solutions}

\begin{equation}\label{m_weak}
\mathcal{M}_{Q}[u,h]=\left\{w\in \left(W_{Q}^{1,2}\right)_0(\Omega):w\geq 0,\;h'(u)w\in \left(W_{Q}^{1,2}\right)_0(\Omega)\right\},
\end{equation}
whenever $h$ is admissible for $u$,
\begin{equation}\label{e_weak}
\mathcal{E}[u]=\left\{\psi^2u:\psi\in C^{1,0}_{0}(\Omega)\right\}
\end{equation}

\begin{equation}\label{a_weak}
\mathcal{A}[u]=\left\{\psi^2h(u)h'(u):\psi\in C^{0,1}_{0}(\Omega),\;h\;\hbox{is admissible for $u$ and}\;h'\geq 0\right\}.
\end{equation}

Basic properties of the above classes as well as the motivation for defining them is summarized below in Remark \ref{weak_classes}. We are now in the position to state the equation that a nonlinear function of a weak (sub-, super-) solution of (\ref{eq}) satisfies (see also \cite{Saw}, \cite{Saw2}).
\begin{Prop}\label{weak_class_prop}
Let $\tilde{u}=h(u)$ and $h$ is admissible for $u$. Then $\tilde{u}$ satisfies
\begin{equation}\label{weak_integral}
-\int(\nabla\tilde{u})^T\tilde{A}\nabla w=-\int (\nabla u)^T\tilde{A}\nabla(wh'(u))+\int w\chi_{\{x\in\Omega:u(x)\notin\mathcal{R}_u\}} h''(u)[\nabla u]^{2}_{\tilde{A}}
\end{equation}
for all $w\in \mathcal{M}_{Q}[u,h]$, where $\tilde{A}(x)=A(x,u(x))$. If in addition $u$ is a weak solution of (\ref{eq}) then $\tilde{u}$ is a $\mathcal{M}_{Q}[u,h]$-weak solution of the following equation
\begin{equation}\label{tilde}
\tilde{L}\tilde{u}=\nabla^T\tilde{A}\nabla\tilde{u}=h'(u)f+\chi_{\{x\in\Omega:u(x)\notin\mathcal{R}_u\}}h''(u)[\nabla u]^{2}_{\tilde{A}}.
\end{equation}
If $u$ is a weak subsolution (supersolution) of (\ref{eq}) and $h'(u)\geq 0(\leq 0)$ then $\tilde{u}$ is a $\mathcal{M}_{Q}[u,h]$-weak subsolution of (\ref{tilde}). Similarly, if $u$ is a weak subsolution (supersolution) of (\ref{eq}) and $h'(u)\leq 0(\geq 0)$ then $\tilde{u}$ is a $\mathcal{M}_{Q}[u,h]$-weak supersolution of (\ref{tilde}).
\end{Prop}

\begin{proof}
Let $u\in W_{Q}^{1,2}$ and $\tilde{u}=h(u)$ with $h$ is admissible for $u$, and let $w\in \mathcal{M}_{Q}[u,h]$. First, we need to verify that all integrals in (\ref{weak_integral}) are absolutely convergent. By Lemma \ref{conv_prop}, or in particular, by Corollary \ref{h_bound} we have that $\nabla\tilde{u}\in \mathcal{L}^2(Q,\Omega)$ and since $h'$ is bounded, $wh'(u)\in \mathcal{L}^2(Q,\Omega)$. Moreover, by the product rule \ref{itm-product-rule} we have that $wh''(u)\nabla_{reg}u\in\mathcal{L}^2(Q,\Omega)$, since according to our convention 
$$
wh''(u)\nabla_{reg}u=\nabla(wh'(u))-(\nabla w)h'(u)\in\mathcal{L}^2(Q,\Omega).
$$
Next, using the chain rule and the product rule we have
\begin{equation*}
\begin{split}
-\int(\nabla\tilde{u})^T\tilde{A}\nabla w&=-\int(\nabla h(u))^T\tilde{A}\nabla w=-\int(h'(u)\nabla u)^T\tilde{A}\nabla w\\
&=-\int (\nabla u)^T\tilde{A}\nabla(wh'(u))+\int w\chi_{\{x\in\Omega:u(x)\notin\mathcal{R}_u\}} h''(u)(\nabla u)^T\tilde{A}\nabla u\\
&=-\int (\nabla u)^T\tilde{A}\nabla(wh'(u))+\int w\chi_{\{x\in\Omega:u(x)\notin\mathcal{R}_u\}} h''(u)[\nabla u]^{2}_{\tilde{A}}
\end{split}
\end{equation*}
The second assertion of the theorem follows immediately.
\end{proof}

The following two lemmas (see \cite{Saw}) give a relationship between different classes of weak solutions.

\begin{Lemma}\label{m_e}
Let $\tilde{u}=h(u)$ and $h$ is admissible for $u$. Then $\mathcal{E}[\tilde{u}]\subset\mathcal{M}_{Q}[u,h]$.
\end{Lemma}

In particular, if $u$ is a weak subsolution (supersolution) of (\ref{eq}) and $h'(u)\geq 0(\leq 0)$ then $\tilde{u}$ is a $\mathcal{E}[\tilde{u}]$-weak subsolution of (\ref{tilde}).

\begin{Lemma}\cite[Lemma 56]{Saw}\label{a_e}
Suppose $u$ is a $\mathcal{A}[u]$-weak subsolution of (\ref{eq}) in $\Omega$, $h$ is admissible for $u$ and $h'(u)\geq 0$. Then $\tilde{u}=h(u)$ is a positive $\mathcal{E}[\tilde{u}]$-weak subsolution of (\ref{tilde}).
\end{Lemma} 

Now, let $w=\psi^2h(u)$ and $\tilde{u}=h(u)$ is an $\mathcal{E}[\tilde{u}]$-weak sub-(super-) solution then from (\ref{tilde}) we obtain

\begin{equation}\label{interm0}
\begin{split}
 - \int \psi^2h(u)h'(u)f \leq (\geq) & \int(\nabla h(u))^T\tilde{A}\nabla \psi^2h(u)\\
 &+\int \chi_{\{x\in\Omega:u(x)\notin\mathcal{R}_u\}}\psi^2h(u) h''(u)[\nabla u]^{2}_{\tilde{A}} .
\end{split}
\end{equation}
For the right-hand side we have

\begin{equation*}
\begin{split}
\int(\nabla h(u))^T\tilde{A}\nabla \psi^2h(u)&+\int \chi_{\{x\in\Omega:u(x)\notin\mathcal{R}_u\}}\psi^2h(u) h''(u)[\nabla u]^{2}_{\tilde{A}}\\
&=\int\psi^2\Gamma(u) [\nabla u]^{2}_{\tilde{A}}+2\int\left\langle \psi\nabla h(u),h(u)\nabla\psi\right\rangle_{\mathcal{L}^2(\Omega,\tilde{A})},
\end{split}
\end{equation*}
where $\Gamma(t)=h'(t)^2+h(t) h''(t)=\left(h(t)^2/2\right)''$. This leads to
\begin{equation}\label{interm2}
\begin{split}
\int\psi^2\Gamma(u) [\nabla u]^{2}_{\tilde{A}}=&\int(\nabla h(u))^T\tilde{A}\nabla \psi^2h(u) -2\int\left\langle \psi\nabla h(u),h(u)\nabla\psi\right\rangle_{\mathcal{L}^2(\Omega,\tilde{A})}\\
& +\int \chi_{\{x\in\Omega:u(x)\notin\mathcal{R}_u\}}\psi^2h(u) h''(u)[\nabla u]^{2}_{\tilde{A}}.
\end{split}
\end{equation}

Applying H\"{o}lder inequality and then Cauchy-Schwarz inequality to the second term, and using (\ref{struc}) we obtain
\begin{equation}\label{interm1}
\begin{split}
2\left|\int\left\langle \psi\nabla h(u),h(u)\nabla\psi\right\rangle_{\mathcal{L}^2(\Omega,\tilde{A})}\right|\leq & \, K\varepsilon \int\psi^2 h'(u)^2[\nabla u]^{2}_{Q}\\
& +K\varepsilon^{-1}\int h(u)^2[\nabla \psi]^{2}_{Q},
\end{split}
\end{equation}
where $0<\varepsilon<1$ . Now let $u\in W_{Q}^{1,2}(\Omega)$ and assume that $\tilde{u}=h(u)$ satisfies one of the following
\begin{itemize}
\item $\Gamma(u)>0$ and $\tilde{u}$ is a $\mathcal{E}[\tilde{u}]$-weak subsolution of (\ref{tilde}) in $\Omega$, or
\item $\Gamma(u)<0$ and $\tilde{u}$ is a $\mathcal{E}[\tilde{u}]$-weak supersolution of (\ref{tilde}) in $\Omega$.
\end{itemize}
Using (\ref{struc}), (\ref{interm0}), (\ref{interm2}) and (\ref{interm1}) we obtain
\begin{equation}\label{e_sub}
\begin{split}
\int\psi^2(|\Gamma(u)|-\varepsilon\frac{K}{k} h'(u)^2)[\nabla u]^{2}_{Q}\leq &\, \varepsilon^{-1}\frac{K}{k}\int h(u)^2[\nabla \psi]^{2}_{Q} \\
&+ \frac{1}{k}\int \psi^2h(u)|h'(u)||f|
\end{split}
\end{equation}

We finish this section with a remark summarizing the properties of the three classes $\mathcal{W}$ of test functions introduced above.
\begin{Rem}\cite[Remark 57]{Saw}\label{weak_classes}
\begin{itemize}
\item If $u$ is a weak solution to (\ref{eq}), then $\mathcal{M}_{Q}[u,h]$ is the maximal subset $\mathcal{W}$ such that $\tilde{u}=h(u)$ is a $\mathcal{W}$-weak solution to (\ref{tilde})
\item The set $\mathcal{E}[\tilde{u}]$ of test functions is the minimal subset $\mathcal{W}$ such that a certain type of Caccioppoli inequality holds for $\mathcal{W}$-weak solutions $\tilde{u}=h(u)$ to (\ref{tilde}) when $u\in W^{1,2}_{Q}(\Omega)$ and $h$ is admissible for $u$ (see inequality (\ref{u_beta}) in Section \ref{Proof of the main result})
\item The set $\mathcal{A}[u]$ of admissible functions is the minimal subset $\mathcal{W}$ such that $\mathcal{W}$-weak solutions to (\ref{eq}) satisfy the property that $\tilde{u}=h(u)$ is a $\mathcal{E}[u]$-weak solution to (\ref{tilde}) for all $h$ that are admissible for $u$.
\end{itemize}
\end{Rem}

\section{Proof of the main result}\label{Proof of the main result}
\subsection{Harnack inequality}
In this section we make use of equation (\ref{tilde}) for nonlinear functions of $u$ and implement the Moser iteration to prove our main technical result, the weak Harnack inequality.
We first prove a Caccioppoli-type inequality for the $Q$-gradient of the powers of $u$. As in \cite{Saw} for $\beta>1$ and $M>0$ we define
\begin{equation*}
h_{\beta,M}(t)=
\begin{cases}
t^\beta, &0<t\leq M,\\
M^\beta+\beta M^{\beta-1}(t-M), &t>M
\end{cases}
\end{equation*}
and $h_{\beta,M}(t)=t^\beta$ for $\beta\leq 1$ and $M>0$. One can show that $h_{\beta,M}$ is admissible for $u\geq m>0$ and then $h_{\beta,M}(u)\in W_{Q}^{1,2}(\Omega)$ since $\Omega$ is bounded (see Corollary \ref{h_bound}). We also have
\begin{equation*}
h_{\beta,M}'(t)=
\begin{cases}
\beta t^{\beta-1}, &0<t\leq M,\\
\beta M^{\beta-1}, &t>M
\end{cases}
\end{equation*}
and 
$$
\Gamma(u)=h_{\beta,M}'(u)^2+h_{\beta,M}(u)h_{\beta,M}''(u)=\eta_{\beta}(u)h_{\beta}'(u)^2
$$
where 
\begin{equation*}
\eta_{\beta}(t)=
\begin{cases}
\frac{2\beta - 1}{\beta}, &0\leq t\leq M\;\hbox{or}\;\beta\leq 1\\
1, &t>M\;\hbox{and}\;\beta> 1
\end{cases}
\end{equation*}

\begin{Lemma}[Caccioppoli inequality]
Let $u\in W_{Q}^{1,2}(\Omega)$ and satisfy one of the following
\begin{subequations}\label{admissible_u}
\begin{align}
&\hbox{$u$ is a positive weak solution of (\ref{eq}) in $\Omega$}\\
&\hbox{$u$ is a positive $\mathcal{A}[u]$-weak subsolution of (\ref{eq}) in $\Omega$ and $\beta>\frac{1}{2}$}\label{admissible_u_b}
\end{align}
\end{subequations}
and $\psi\in C_{0}^{0,1}(\Omega)$, $\beta\in\RR$ with $\beta\neq 0,\frac{1}{2}$. Then the following inequality holds
\begin{equation}\label{u_beta}
\int\psi^2[\nabla u^{\beta}]_{Q}^2\leq C\left(\mu_{\beta}^{-2}\int u^{2\beta}[\nabla\psi]^2+\mu_{\beta}^{-1}|\beta|\int\psi^2u^{2\beta-1}|f|\right)
\end{equation}
where
\begin{equation}\label{mu}
\mu_{\beta}=\min\left\{\left|\frac{2\beta-1}{\beta}\right|,1\right\}.
\end{equation}
\end{Lemma}

\begin{proof}
It follows from Proposition \ref{weak_class_prop} and Lemma \ref{m_e} that $\tilde{u}=h(u)=h_{\beta,M}(u)$ is a positive $\mathcal{E}[\tilde{u}]$-weak solution of (\ref{tilde}) provided (\ref{admissible_u}a) holds. On the other hand, by Lemma \ref{a_e} $\tilde{u}$ is a positive $\mathcal{E}[\tilde{u}]$-weak subsolution of (\ref{tilde}) if (\ref{admissible_u}a) holds (and therefore $\Gamma>0$). Therefore, inequality (\ref{e_sub}) holds.
Next, we note that $|\Gamma(u)|\geq \mu_{\beta}h'(u)^2$ where $\mu_{\beta}$ is defined by (\ref{mu}) and choose $\varepsilon=\frac{1}{2}\frac{k}{K}\mu_\beta$ to obtain
$$
\frac{1}{2}\mu_\beta\int\psi^2 h'(u)^2[\nabla u]_{Q}^2 \leq C\mu_{\beta}^{-1}\int h(u)^2[\nabla \psi]_{Q}^2 + C\int\psi^2 h(u)|h'(u)||f|
$$
Letting $M\rightarrow\infty$ yields the above estimate with $u^{\beta}$ in place of $h_{\beta,M}(u)$. Finally, using $\left((u^{\beta})'\right)^2[\nabla u]_{Q}^2=[\nabla u^{\beta}]_{Q}^2$ we have
$$
\int\psi^2[\nabla u^{\beta}]_{Q}^2\leq C\left(\mu_{\beta}^{-2}\int u^{2\beta}[\nabla\psi]_{Q}^2+\mu_{\beta}^{-1}|\beta|\int\psi^2u^{2\beta-1}|f|\right)
$$
whenever $\psi\in C_{0}^{0,1}(\Omega)$, $\beta\in\RR$ with $\beta\neq 0,\frac{1}{2}$, $u$ is bounded below by a positive constant, and it satisfies (\ref{admissible_u}).
\end{proof}

\begin{Rem}
The assumption (\ref{admissible_u_b}) and the $\mathcal{M}_Q[u,h]$ class of weak solutions are only used in the proof of the local boundedness result, Proposition \ref{local_boundedness}.
\end{Rem}


We now perform the Moser iteration on a weak solution of (\ref{eq}), using the accumulating sequence of cutoff functions (\ref{cutoff}) to obtain a local bound on $\sup{\tilde{u}}$ by an $L^2$-norm of $\tilde{u}$. This is a crucial step in proving Harnack inequality.

\begin{Prop}[Moser iteration]\label{Moser iteration}
Let $|\gamma|\leq 2$, $u\in W_{Q}^{1,2}(\Omega)$ satisfy one of the following conditions
\begin{subequations}\label{admissible_u2}
\begin{align}
&\hbox{$u$ is a nonnegative weak solution of (\ref{eq}) in $\Omega$}\\
&\hbox{$u+m$ is a positive $\mathcal{A}[u+m]$-weak subsolution of (\ref{eq}) in $\Omega$ for all $m>0$ and $\gamma>\frac{1}{2}$}.
\end{align}
\end{subequations}
Then for any metric ball $B=B(y,r)$ with $y\in\Omega$, $0<r<\eta\, dist(y,\partial\Omega)$ the following inequality holds
\begin{equation}\label{pre_harnack}
\esssup_{x\in B(y,\nu r)}\overline{u}(x)^\gamma\leq \frac{C_\sigma}{(1-\nu)^{\frac{1}{(\sigma-1)^2}+1}(\delta_y(\nu r)/r)^{\frac{\sigma}{\sigma-1}}} \left\{\frac{1}{|B|}\int_{B}\overline{u}^{2\gamma}\right\}^{\frac{1}{2}}
\end{equation}
where $\overline{u}=u+m(r)$ with $m(r)=r^2||f||_{L^{\infty}}$, $\tau=\sigma/(\sigma-1)$ with $\sigma$ and $\nu$ as in (\ref{sobolev}) and (\ref{cutoff}). The constant $C_{\sigma}$ depends only on $\sigma$ and the constant $C$ from (\ref{cutoff}).
\end{Prop}

\begin{proof}
First assume that $\overline{u}$ satisfies (\ref{admissible_u}), we will later see that for our choice of $\beta$ it follows from (\ref{admissible_u2}). Therefore, we can apply (\ref{u_beta}) and for each $\psi\in C_{0}^{0,1}(\Omega)$ and $\beta\in\RR$ with $\beta\neq 0,\frac{1}{2}$ we have
\begin{equation*}
\begin{split}
\int\psi^2[\nabla \overline{u}^{\beta}]_{Q}^2&\leq C\left(\mu_{\beta}^{-2}\int \overline{u}^{2\beta}[\nabla\psi]_{Q}^2+\mu_{\beta}^{-1}|\beta|\int\psi^2\overline{u}^{2\beta-1}|f|\right)\\
&\leq C\left(\mu_{\beta}^{-2}\int \overline{u}^{2\beta}[\nabla\psi]_{Q}^2+\mu_{\beta}^{-1}|\beta|\int\psi^2\overline{u}^{2\beta}m(r)^{-1}|f|\right),
\end{split}
\end{equation*}
where the last inequality is due to the fact that $\overline{u}^{-1}\leq m(r)^{-1}$. Noting that $m(r)^{-1}|f|\leq r^{-2}$ we get
\begin{equation*}
\int\psi^2[\nabla \overline{u}^{\beta}]_{Q}^2 \leq C\left\{\mu_{\beta}^{-2}\int \overline{u}^{2\beta}[\nabla\psi]_{Q}^2+\mu_{\beta}^{-1}|\beta|r^{-2}\int\psi^2\overline{u}^{2\beta}\right\}.
\end{equation*}
Next, we use the above inequality with an accumulating sequence of cutoff functions $\psi_j$ and a fixed subunit ball $B=B(y,r)$. Recall that $E_j=\supp(\psi_j) = B(y,r_j)$ and $\overline{u}>0$ on $E_j$ $\forall j$. Dividing through by the measure of the support of $\psi_j$, $|E_j|$, and taking square roots we obtain
\begin{equation}\label{ineq1}
\begin{array}{rcl}
\displaystyle{\left\{\frac{1}{|E_j|}\int_{E_j}\psi_j^2[\nabla \overline{u}^{\beta}]_{Q}^2\right\}^{\frac{1}{2}}} & \leq & \displaystyle{C\left\{\mu_{\beta}^{-2}\frac{1}{|E_j|}\int_{E_j} \overline{u}^{2\beta}[\nabla\psi_j]_{Q}^2\right\}^{\frac{1}{2}} }\\
&&+\displaystyle{C\left\{\mu_{\beta}^{-1}|\beta|r^{-2}\frac{1}{|E_j|}\int_{E_j}\psi_j^2\overline{u}^{2\beta}\right\}^{\frac{1}{2}}.}
\end{array}
\end{equation}
We now use the Sobolev inequality (\ref{sobolev}) for $w=\psi_j\overline{u}^{\beta}$. Then, for some $\sigma>1$ we obtain using (\ref{ineq1})
\begin{equation*}
\begin{split}
\left\{\frac{1}{|E_j|}\int_{E_j}(\psi_j\overline{u}^\beta)^{2\sigma}\right\}^{\frac{1}{2\sigma}}&\leq
Cr_j\left\{\frac{1}{|E_j|}\int_{E_j}[\nabla(\psi_j\overline{u}^\beta)]_{Q}^2\right\}^{\frac{1}{2}}
+C\left\{\frac{1}{|E_j|}\int_{E_j}(\psi_j\overline{u}^\beta)^2\right\}^{\frac{1}{2}}\\
&\leq Cr\left(1+\mu_{\beta}^{-1}\right)\left\{\frac{1}{|E_j|}\int_{E_j}\overline{u}^{2\beta}[\nabla\psi_j]_{Q}^2\right\}^{\frac{1}{2}}\\
&\quad+C\left(1+\mu_{\beta}^{-1/2}|\beta|^{1/2}\right)\left\{\frac{1}{|E_j|}\int_{E_j}(\overline{u}^{\beta}\psi_j)^2\right\}^{\frac{1}{2}}.
\end{split}
\end{equation*}
Using the property of the cutoff function $\psi_j$
$$
[\nabla\psi_j]_{Q}\leq \frac{C}{(1-\nu)\delta(\nu r)\left(1-\frac{\delta_y(\nu r)}{r}\right)^{j}}
$$
and writing $\delta(\nu r)$ for $\delta_y(\nu r)$ yields
\begin{equation}\label{ineq_to_iter}
\begin{array}{l}
\displaystyle{\left\{\frac{1}{|E_j|}\int_{E_{j+1}}(\overline{u}^\beta)^{2\sigma}\right\}^{\frac{1}{2\sigma}} }\\
\leq \displaystyle{C\left(\frac{(1+\mu_{\beta}^{-1})r}{(1-\nu)\delta(\nu r)\left(1-\frac{\delta(\nu r)}{r}\right)^{j}}
+\mu_{\beta}^{-1/2}|\beta|^{1/2}+1\right)\left\{\frac{1}{|E_j|}\int_{E_j}\overline{u}^{2\beta}\right\}^{\frac{1}{2}}.}
\end{array}
\end{equation}
It can be shown that the constant in the parentheses is bounded by
$$
\frac{r}{(1-\nu)\delta(\nu r)\left(1-\frac{\delta(\nu r)}{r}\right)^{j}}\cdot M(\beta)
$$
with
$$
M(\beta):=C\left(1+|\beta|^{1/2}+\frac{|\beta|}{|2\beta-1|^{1/2}}+\frac{|\beta|}{|2\beta-1|}\right).
$$
Next, to iterate inequality (\ref{ineq_to_iter}), we fix $\gamma\neq 0$ and denote $u_j=\overline{u}^{\gamma\sigma^{j-1}}$. Take $\beta=\gamma\sigma^{j-1}$, one can see that for each $j\geq 1$ (\ref{admissible_u}) follows from (\ref{admissible_u2}) for this choice of $\beta$. Using $u_{j+1}=u_j^\sigma$ we obtain from (\ref{ineq_to_iter})
\begin{equation}\label{ineq2}
\begin{array}{l}
\displaystyle{\left\{\frac{1}{|E_{j+1}|}\int_{E_{j+1}}u_{j+1}^2\right\}^{\frac{1}{2}}} \\
\leq\displaystyle{\frac{|E_j|^{1/2}}{|E_{j+1}|^{1/2}}C\left(\frac{M(\gamma\sigma^{j-1})r}{(1-\nu)\delta(\nu r)\left(1-\frac{\delta(\nu r)}{r}\right)^{j}}\right)^\sigma\left\{\frac{1}{|E_j|}\int_{E_j}u_j^{2}\right\}^{\frac{\sigma}{2}}.}
\end{array}
\end{equation}

Let
$$
N_j=\left\{\frac{1}{|E_{j}|}\int_{E_{j}}u_{j}^2\right\}^{\frac{1}{2\sigma^{j-1}}},
$$

then from (\ref{ineq2})

\begin{equation*}
N_{j+1}\leq\left[\frac{|E_j|^{1/2}}{|E_{j+1}|^{1/2}}C\left(\frac{M(\gamma\sigma^{j-1})r}{(1-\nu)\delta(\nu r)\left(1-\frac{\delta(\nu r)}{r}\right)^{j}}\right)^\sigma\right]^{\frac{1}{\sigma^j}}N_j
\end{equation*}

Iterating, we obtain

$$
\limsup_{j\rightarrow\infty}N_j\leq
C\left\{\prod\limits_{j=1}^{\infty}\left[\frac{|E_j|^{1/2}}{|E_{j+1}|^{1/2}}C\left(\frac{M(\gamma\sigma^{j-1})r}{(1-\nu)^{j+1}\delta(\nu r)}\right)^\sigma\right]^{\frac{1}{\sigma^j}}\right\}N_1,
$$
where we also used the fact that $1-\delta(\nu r)/r\geq 1-\nu$.
Next, recall the properties of the supports of $\psi_j$
$$
\frac{|E_j|}{|E_{j+1}|}\leq C.
$$

Therefore

\begin{equation}\label{iter_result1}
\limsup_{j\rightarrow\infty}N_j\leq C\frac{r^{\frac{\sigma}{\sigma-1}}}{(1-\nu)^{\frac{1}{(\sigma-1)^2}+1}\delta(\nu r)^{\frac{\sigma}{\sigma-1}}}
\prod\limits_{j=1}^{\infty}(M(\gamma\sigma^{j-1}))^{\frac{1}{\sigma^{j-1}}}\left\{\frac{1}{|E_1|}\int_{E_1}\overline{u}^{2\gamma}\right\}^{\frac{1}{2}}
\end{equation}
where we used
\begin{equation*}
\prod\limits_{j=1}^{\infty}\frac{1}{\left(1-\nu\right)^{j/\sigma^{j-1}}}
=exp\left(\sum\limits_{j=1}^{\infty}-\frac{j}{\sigma^{j-1}}\ln \left(1-\nu\right)\right)
=(1-\nu)^{-\frac{1}{(\sigma-1)^2}}.
\end{equation*}
Also 

\begin{equation}\label{iter_result2}
\esssup_{x\in B(y,\nu r)}\overline{u}^\gamma \leq \limsup_{j\rightarrow\infty}N_j
\end{equation}

and therefore we only need to estimate the constant 

$$
C_\gamma=\prod\limits_{j=1}^{\infty}\left(M(\gamma\sigma^{j-1})\right)^{\frac{1}{\sigma^{j-1}}},
$$

where 

$$
M(\gamma\sigma^{j-1})= C\left(1+|\gamma\sigma^{j-1}|^{1/2}+\frac{|\gamma\sigma^{j-1}|}{|2\gamma\sigma^{j-1}-1|^{1/2}}+\frac{|\gamma\sigma^{j-1}|}{|2\gamma\sigma^{j-1}-1|}\right).
$$

We can write

$$
C_\gamma=C\exp\sum\limits_{j=1}^{\infty}\frac{\ln \left[j^N\left(1+|\gamma\sigma^{j-1}|^{1/2}+\frac{|\gamma\sigma^{j-1}|}{|2\gamma\sigma^{j-1}-1|^{1/2}}+\frac{|\gamma\sigma^{j-1}|}{|2\gamma\sigma^{j-1}-1|}\right)\right]}{\sigma^{j-1}}.
$$

In order to have $C_\gamma\leq\infty$ we must ensure that $\gamma\sigma^{j-1}$ does not approach $1/2$ for any values of $j$. We first note that for any $\gamma>0$ we have $\tilde{\gamma}\leq\gamma<\tilde{\gamma}\sigma$ with $\tilde{\gamma}$ of the form $\tilde{\gamma}=\frac{1}{4}\sigma^k(\sigma+1)$ for some integer $k$. Therefore due to the monotonicity of 
$$
\left(\frac{1}{|B|}\int_{B}\overline{u}^\gamma\right)^{1/\gamma}
$$
it is enough to prove (\ref{pre_harnack}) for $\gamma=\tilde{\gamma}$. For this choice of $\gamma$ we have
$$
|2\tilde{\gamma}\sigma^{j-1}-1|\geq\frac{1}{2}(1-\sigma^{-1})
$$

Indeed, if $k\geq -j+1$ then
$$
2\tilde{\gamma}\sigma^{j-1}-1=\frac{1}{2}\sigma^{k+j-1}(\sigma+1)-1\geq\frac{\sigma}{2}-\frac{1}{2}\geq\frac{1}{2}(1-\sigma^{-1})
$$
while for $k\leq -j$ we similarly check that
$$
2\tilde{\gamma}\sigma^{j-1}-1\leq -\frac{1}{2}(1-\sigma^{-1})
$$

We therefore obtain

$$
C_\gamma\leq C_\sigma(1+|\gamma|^{p(\sigma)})
$$

and for $|\gamma|\leq 2$ from (\ref{iter_result1}), (\ref{iter_result2})

\begin{equation}\label{moser_result}
\begin{split}
\esssup_{x\in B(y,\nu r)}\overline{u}^\gamma &\leq \frac{C_\sigma}{(1-\nu)^{\frac{1}{(\sigma-1)^2}+1}(\delta(\nu r)/r)^{\frac{\sigma}{\sigma-1}}} \left\{\frac{1}{|E_1|}\int_{E_1}\overline{u}^{2\gamma}\right\}^{\frac{1}{2}}\\
&=\frac{C_\sigma}{(1-\nu)^{\frac{1}{(\sigma-1)^2}+1}(\delta(\nu r)/r)^{\frac{\sigma}{\sigma-1}}} \left\{\frac{1}{|B|}\int_{B}\overline{u}^{2\gamma}\right\}^{\frac{1}{2}}.
\end{split}
\end{equation}
Finally, we note that if $\gamma>1/2$ then for $\beta=\gamma\sigma^{j-1}$ we also have $\beta>1/2$. Therefore, $\overline{u}=u+m$ satisfies (\ref{admissible_u}) as long as $u$ satisfies (\ref{admissible_u2}).
\end{proof}

To establish a weak Harnack inequality we will adopt an argument of Bombieri (see \cite[Lemma 3]{Mos}). One needs to be careful however, since in this case the coefficients in the inequalities depend on the radius of the ball.

\begin{Lemma}\label{bound}
Let $w>0$ be a measurable function defined in a neighborhood of $B_0=B(y_0,r_0)$. Suppose there exist positive constants $\tau$ and $\nu_0<1$, and $a\geq 0$; and decreasing functions $c_1(r)$, $c_2(r)$, with $0<c_1(r),c_2(r)<\infty$ for any $0<r<\infty$, such that the following two conditions hold
\begin{enumerate}
\item 
\begin{equation}\label{moser}
\sup_{x\in \nu B}w^{\gamma}\leq \frac{c_1(\nu r)}{(1-\nu)^{\tau}}\left\{\frac{1}{|B|}\int\limits_{B}w^{2\gamma}\right\}^{1/2}
\end{equation}
for every $0<\nu_0<\nu<1$, $0<\gamma\leq 2$, where $B=B(y,r)\subset B_0$, and 
\item 
\begin{equation}\label{weak}
s|\{x\in B: \log w > s+a\}|<c_2(r)|B|
\end{equation}
for every $s>0$.
\end{enumerate}
Then there exists $b=b(\nu_0,c_1(r),c_2(r))$ such that
\begin{equation}\label{exp_bound}
\esssup_{B(y,\nu_0 r)}w<b e^a.
\end{equation}
\end{Lemma}

\begin{proof}
Define
$$
\varphi(\rho)=\esssup_{y\in B(x,\rho)}(\log w(y)-a)\quad\hbox{for}\quad \nu_0r\leq\rho\leq r.
$$
First, note that if $\varphi(\nu_0 r)\leq 0$ then estimate (\ref{exp_bound}) holds with $b\equiv const$, therefore, we may assume $\varphi(\rho)>0$ for all $\nu_0r\leq\rho\leq r$. We can decompose the ball $B$ in the following way
$$
B=\{y\in B: \log w(y)-a>\frac{1}{2}\varphi(r)\}\cup\{y\in B: \log w(y)-a\leq\frac{1}{2}\varphi(r)\}.
$$
Then, using condition (\ref{weak}) we have
$$
e^{-2\gamma a}\int\limits_{B}w^{2\gamma}=\int\limits_{B}e^{2\gamma(\log w-a)}\leq \frac{2c_2(r)}{\varphi}e^{2\gamma\varphi}|B| +e^{\gamma\varphi}|B|
$$ 
or 
$$
\frac{e^{-2\gamma a}}{| B|}\int\limits_{B}w^{2\gamma}\leq \frac{2c_2(r)}{\varphi}e^{2\gamma\varphi} +e^{\gamma\varphi}
$$ 
where $\varphi = \varphi(r)$.

Next, choose $\gamma$ so that the two terms on the right-hand side are equal ${(2c_2(r)/\varphi)e^{2\gamma\varphi}=e^{\gamma\varphi}}$, i.e.
\begin{equation}\label{gamma}
\gamma = \frac{1}{\varphi}\log\left(\frac{\varphi}{2c_2(r)}\right)
\end{equation}
To satisfy the condition $0<\gamma\leq 2$ we must then require 
\begin{equation}\label{req1}
\varphi = \varphi(r)>c_3 c_2(r).
\end{equation}
In that case, we have
$$
\left(\frac{e^{-2\gamma a}}{|B|}\int\limits_{B}w^{2\gamma}\right)^{\frac{1}{2}}\leq \sqrt{2} e^{\gamma\varphi/2}
$$
and using (\ref{moser})
$$
\varphi(\nu r)\leq\frac{1}{\gamma}\log\left(\frac{c_1(\nu r)\sqrt{2}}{(1-\nu)^{\tau}}e^{\gamma\varphi/2}\right)=\frac{1}{\gamma}\log\left(\frac{c_1(\nu r)\sqrt{2}}{(1-\nu)^{\tau}}\right)+\frac{1}{2}\varphi(r)
$$
hence, by (\ref{gamma}),
$$
\varphi(\nu r)<\varphi(r)\left(\frac{\log(c_1(\nu r)\sqrt{2}/(1-\nu)^{\tau})}{\log(\varphi/2c_{2}(r))}+\frac{1}{2}\right).
$$
If 
\begin{equation}\label{req2}
\varphi(r)>\frac{8c_{1}(\nu r)^{4}c_2(r)}{(1-\nu)^{4\tau}}
\end{equation}
then the first term in the parentheses is less than $1/4$ and therefore we have $\varphi(\nu r)<3/4\varphi(r)$. On the other hand, if (\ref{req1}) or (\ref{req2}) do not hold, then there exists $c_4$ such that
$$
\varphi(r)\leq\frac{c_4c_{1}(\nu r)^{4}c_2(r)}{(1-\nu)^{4\tau}}.
$$
Therefore, in any case we have
\begin{equation}\label{ineq}
\varphi(\nu r)<\frac{3}{4}\varphi(r)+\frac{c_4c_{1}(\nu r)^{4}c_2(r)}{(1-\nu)^{4p\tau}}
\end{equation}
since $\varphi(\nu r)\leq\varphi(r)$. We now iterate this inequality with 
$$
0<\nu_0<\nu_1<\nu_2<\ldots<\nu_k\leq 1
$$ 
to obtain
\begin{equation*}
\begin{split}
\varphi(\nu_0r)&<\left(\frac{3}{4}\right)^k\varphi(r)+c_4\sum\limits_{j=0}^{k-1}\left(\frac{3}{4}\right)^jc_{1}(\nu_{j}r)^{4}c_2(\nu_{j+1}r)(\nu_{j+1}-\nu_j)^{-4\tau}\\
&\leq\left(\frac{3}{4}\right)^k\varphi(r)+c_4c_{1}(\nu_{0}r)^{4}c_2(\nu_{0}r)\sum\limits_{j=0}^{k-1}\left(\frac{3}{4}\right)^j(\nu_{j+1}-\nu_j)^{-4\tau}
\end{split}
\end{equation*}
where the last inequality is due to the fact that $c_1(r)$ and $c_2(r)$ are decreasing.
Choosing 
$$
\nu_j=1-\frac{1-\nu_0}{1+j}
$$
and letting $k\rightarrow\infty$ we obtain

$$
\varphi(\nu_0r)\leq C(\nu_0,\tau)c_{1}(\nu_{0}r)^{4}c_2(\nu_{0}r)
$$
and therefore
$$
\esssup_{B(y,\nu_0 r)}w<C(\nu_0,\tau)\exp(c_{1}(\nu_{0}r)^{4}c_2(\nu_{0}r)) e^a
$$
which implies (\ref{exp_bound}) with $b=C(\nu_0,\tau)\exp(c_{1}(\nu_{0}r)^{4}c_2(\nu_{0}r))$.
\end{proof}

We now want to apply Lemma \ref{bound} to a positive subsolution $\overline{u}$. It has been shown in Proposition \ref{Moser iteration} that condition (\ref{moser}) is satisfied with $c_1(\nu r)=C(\nu r/\delta(\nu r))^{\sigma/(\sigma-1)}$ and $\tau=1/(\sigma-1)^2+1$. The following lemma shows that (\ref{weak}) holds for $\overline{u}$ and $1/\overline{u}$. 

\begin{Lemma}\label{lemma_weak1}
Let $u\in W_{Q}^{1,2}(\Omega)$ be a nonnegative weak solution of (\ref{eq}) in $\Omega$, and let $\overline{u}=u+m(r)$, $m(r)=r^2||f||_{L^{\infty}}$ with $r\leq r_0$. There holds
\begin{equation}\label{weak_est1}
s|\{x\in B: \log \overline{u} > s+\left\langle \log \overline{u}\right\rangle_{B}\}|<C\frac{|B|r}{\delta(r)}, and
\end{equation}
\begin{equation}\label{weak_est2}
s|\{x\in B: \log(1/ \overline{u}) > s-\left\langle \log \overline{u}\right\rangle_{B}\}|<C\frac{|B|r}{\delta(r)}
\end{equation}
where $B=B(y,r)\subset B_0$, $B(y,r+\delta(r))\subset B_0$, $y\in\Omega$, $0<r<\eta\, dist(y,\partial\Omega)$, and $\delta(r)$ as in (\ref{nondoub_order}).
\end{Lemma}

\begin{proof}
Denote $v=\ln\overline{u}$, for any $s>0$, it is easy to check that $v\in  W_{Q}^{1,2}(\Omega)$, and we have
$$
s|\{x\in B: v-\left\langle v\right\rangle_B > s\}|\leq\int_{B}|v-\left\langle v\right\rangle_B|.
$$
Applying the Poincar\'{e} inequality (\ref{poincare}) we obtain
$$
s|\{x\in B: v-\left\langle v\right\rangle_B > s\}|\leq Cr \int\limits_{B}[\nabla v]_{Q}.
$$
Therefore, in order to prove (\ref{weak_est1}) it is enough to show
\begin{equation}\label{weak_est_inter}
\int\limits_{B}[\nabla v]_{Q}\leq \frac{C}{\delta(r)}|B|.
\end{equation}
Consider equation (\ref{sol}) and substitute $w=\frac{\varphi^2}{\overline{u}}$ with $\varphi\in W^{1,2}_{0}(B(y,r+\delta(r)))$ as in (\ref{spec_cutoff}). Using the third property in (\ref{spec_cutoff}) we obtain 
\begin{equation*}
\begin{split}
\int_{B(y,r+\delta(r))}\varphi^2(\nabla\ln(\overline{u}))^T\tilde{A}\nabla\ln(\overline{u})\\
=2\int_{B(y,r+\delta(r))}\varphi\nabla\varphi \tilde{A}\nabla\ln(\overline{u})+\int_{B(y,r+\delta(r))}\frac{f\varphi^2}{\overline{u}}\\
\leq \frac{C}{\delta(r)}\int_{B(y,r+\delta(r))}\varphi[\nabla\ln(\overline{u})]_{\tilde{A}}+\int_{B(y,r+\delta(r))}\frac{\varphi^2}{r^2}
\end{split}
\end{equation*}
Using H\"{o}lder inequality and Cauchy-Schwarz inequality we have
$$
\int_{B(y,r+\delta(r))}\varphi[\ln(\overline{u})]_{\tilde{A}}\leq \frac{C}{\varepsilon \delta(r)}|B(y,r+\delta(r))| + \varepsilon \delta(r)\int_{B(y,r+\delta(r))}\varphi^2(\nabla\ln(\overline{u}))^T\tilde{A}\nabla\ln(\overline{u})
$$
and using the above we get
\begin{align*}
\int_{B(y,r+\delta(r))}\varphi[\ln(\overline{u})]_{\tilde{A}}\leq &\frac{C}{\varepsilon \delta(r)}|B(y,r+\delta(r))|\\
&+C\varepsilon\int_{B(y,r+\delta(r))}\varphi[\nabla\ln(\overline{u})]_{\tilde{A}}+\frac{\varepsilon \delta(r)}{r^2}|B(y,r+\delta(r))|
\end{align*}
Now, we choose $\varepsilon$ small enough that the second term on the right can be absorbed into the left-hand side and use (\ref{struc}) and the fact that $\delta(r)<r$ to obtain
$$
\int\limits_{B}[\nabla v]_{Q}\leq \int_{B(y,r+\delta(r))}\varphi[\ln(\overline{u})]_{Q}
\leq K\int_{B(y,r+\delta(r))}\varphi[\ln(\overline{u})]_{\tilde{A}}\leq \frac{C}{\delta(r)}|B(y,r+\delta(r))|
$$
By the definition of $\delta(r)$, (\ref{nondoub_order}), this implies (\ref{weak_est_inter}).\\
The proof of (\ref{weak_est2}) follows in a similar way.
\end{proof}

We can now establish a weak version of Harnack inequality

\begin{Th}
Let $u$ be a nonnegative weak solution of (\ref{eq}) in $B(y,r)$. Then $u$ satisfies the following Harnack inequality
\begin{equation}\label{harnack}
\esssup_{x\in B(y,\nu_0r)}(u(x)+m(r))\leq C_{Har}(r) \essinf_{x\in B(y,\nu_0r)}(u(x)+m(r))
\end{equation}
where $C_{Har}(r)=C\exp\left(2\left[\frac{\nu_0r}{\delta(\nu_0 r)}\right]^{\frac{4\sigma}{\sigma-1}+1}\right)$, so $C_{Har}\to\infty$ as $r\to 0$ provided $\delta(r)/r\to 0$ as $r\to 0$, and $m(r)=r^2||f||_{L^{\infty}}$.
\end{Th}
\begin{proof}
It follows from (\ref{pre_harnack}), (\ref{weak_est1}), and (\ref{weak_est2}) that the functions $w_1=\overline{u}$ and $w_2=1/\overline{u}$ satisfy the conditions of Lemma \ref{bound} with $a_1=\left\langle\ln\overline{u}\right\rangle_B$ and $a_2=-\left\langle\ln\overline{u}\right\rangle_B$ respectively, and
$$
c_1(\nu r)=\left[\frac{\nu r}{\delta(\nu r)}\right]^{\frac{\sigma}{\sigma-1}}
$$
$$
c_2(r)=\frac{r}{\delta(r)}
$$
Therefore we have
$$
\esssup_{\nu_0B}\overline{u}<C e^{\left\langle\ln\overline{u}\right\rangle_B}\exp(c_{1}(\nu_{0}r)^{4}c_2(\nu_{0}r))
$$
$$
\esssup_{\nu_0B}\frac{1}{\overline{u}}<C e^{-\left\langle\ln\overline{u}\right\rangle_B}\exp(c_{1}(\nu_{0}r)^{4}c_2(\nu_{0}r))
$$
Multiplying the above inequalities we obtain (\ref{harnack}).
\end{proof}

\subsection{Local boundedness of weak solutions}
In this section we prove local boundedness of weak solutions. Note that the assumptions needed to establish local boundedness are weaker then the ones we require for the continuity result. In particular, we do not require any version of the Poincar\'{e} inequality, only the weak Sobolev inequality (\ref{sobolev}). The proof is similar to the proof of Proposition 59 in \cite{Saw} and will follow from (\ref{pre_harnack}).
\begin{Prop}[Local boundedness]\label{local_boundedness}
Suppose that $u$ is a weak solution to (\ref{eq}) in $\Omega$ and that (\ref{contain}), (\ref{sobolev}), and (\ref{cutoff}) hold. Then $u$ is locally bounded in $\Omega$. More precisely, if $y\in\Omega$, $0<r<\eta dist(y,\partial\Omega)$ for sufficiently small $\eta$, then
\begin{equation}\label{boundedness}
||u||_{L^{\infty}(B(y,\nu r))}\leq C(r)\left\{\left(\frac{1}{|B(y, r)|}\int_{B(y,r)}u^2\right)^{\frac{1}{2}}+r^2||f||_{L^\infty}\right\}
\end{equation}
with $C(r)=C(\delta(\nu r)/r)^{\frac{\sigma}{\sigma-1}}<\infty$ $\forall r>0$, and $\nu>0$ is as in (\ref{pre_harnack}).
\end{Prop}
The proof follows from the proof of \cite[Proposition 59]{Saw}, but the argument is simpler in our case and we provide it here for the sake of completeness.
\begin{proof}
Let $m(r)=r^2||f||_{\infty}$ if $||f||_{\infty}\neq 0$ and $m(r)=m>0$ if $||f||_{\infty}= 0$. In the latter case we will take the limit $m\to 0$ in the end of the proof. Consider $h(t)=\sqrt{m(r)^2+t^2}$ for $t\geq 0$ and $h(t)=m(r)$ for $t<0$. It is easy to check that $h(t)$ is admissible for $u$ and $h'(t)\geq 0$. Therefore, by Proposition \ref{weak_class_prop} $\tilde{u}=h(u)$ is a positive $\mathcal{M}_{Q}[u,h]$-weak subsolution to (\ref{tilde}):
\begin{equation*}
\tilde{L}\tilde{u}=\nabla^T\tilde{A}\nabla\tilde{u}=h'(u)f+\chi_{\{x\in\Omega:u(x)\notin\mathcal{R}_u\}}h''(u)[\nabla u]^{2}_{\tilde{A}}
\end{equation*}
where
\begin{equation*}
\mathcal{M}_{Q}[u,h]=\left\{w\in \left(W_{Q}^{1,2}\right)_0(\Omega):w\geq 0,\;h'(u)w\in \left(W_{Q}^{1,2}\right)_0(\Omega)\right\}
\end{equation*}
Moreover, since $h''(t)\geq 0$ we have that $\tilde{u}$ is a positive $\mathcal{M}_{Q}[u,h]$-weak subsolution to the following equation
\begin{equation}\label{m_eq}
\tilde{L}\tilde{u}=\nabla^T\tilde{A}\nabla\tilde{u}=h'(u)f
\end{equation}
Next, we want to apply inequality (\ref{pre_harnack}) to $\tilde{u}+m$. In order to do this we need to show that $\tilde{u}+m$ is a positive $\mathcal{A}[\tilde{u}+m]$-weak subsolution of (\ref{m_eq}). We have already shown that $\tilde{u}$ (and therefore $\tilde{u}+m$) is a positive $\mathcal{M}_{Q}[u,h]$-weak subsolution, therefore it is enough to show that $\mathcal{A}[\tilde{u}+m]\subset\mathcal{M}_{Q}[u,h]$. This can be done exactly as in the proof of \cite[Proposition 59]{Saw} and we omit it here. Applying inequality (\ref{pre_harnack}) to $\tilde{u}+m$ we then have
\begin{equation*}
||u||_{L^\infty(B(y,\nu r))}\leq ||\overline{\tilde{u}}||_{L^\infty(B(y,\nu r))}\leq C(r)^{\frac{1}{\gamma}}_{\nu, \tau}\left\{\frac{1}{|B(y,r)|}\int_{B(y,r)}\overline{\tilde{u}}^{2\gamma}\right\}^{\frac{1}{2\gamma}}
\end{equation*}
for all $\gamma>1/2$ and, in particular, for $\gamma=1$. Here we denoted $\overline{\tilde{u}}=\tilde{u}+m(r)$. To finish the proof we only need to observe, that $\overline{\tilde{u}}^2=u^2+m(r)^2+2m(r)\tilde{u}+m(r)^2\leq 2u^2 +3m(r)^2$ which concludes (\ref{boundedness}).
\end{proof}
A similar result has been recently established in a more general setting \cite{Wheeden}.

\subsection{Continuity of weak solutions}

We will now establish continuity of solutions under the condition of a certain control on the non doubling order of metric balls $\delta(r)$.
The next proposition follows from Harnack inequality and is a weak version of \cite[Theorem 8.22]{Gilb}, since the constant in Harnack inequality is not uniform in radius.
\begin{Prop}
Let $u$ be a nonnegative weak solution of (\ref{eq}) in $B(y,R)$. For $0<r<R$ denote 
$$
\omega(r)=\esssup_{x\in B(y,r)}u(x)-\essinf_{x\in B(y,r)}u(x)
$$
Then we have 
\begin{equation}\label{pre_holder}
\omega(r)\leq C\omega(R)\left(\frac{r}{R}\right)^{\alpha(r)}+C\,C_{Har}(r)\left(\frac{r}{R}\right)^{2\mu}, \quad 0<r\leq R
\end{equation}
where $\mu\in(0,1)$, $\alpha(r)\to 0$ as $r\to 0$, and $C_{Har}(r)$ as in (\ref{harnack}).
\end{Prop}

\begin{proof}
We first claim, that 
\begin{equation}\label{pre_continuity}
\omega(\nu_0 r)\leq\left(1-\frac{1}{2C_{Har}(r)}\right)\omega(r)+r^2||f||_{\infty}
\end{equation}
This follows in exactly the same way as in \cite[Section 3.4]{Saw} and we therefore omit the proof. Now, let $\gamma(r)=1-1/(2C_{Har}(r))$ and we proceed to iterate inequality \ref{pre_continuity} in a similar way as in the proof of \cite[Lemma 8.23]{Gilb}. Namely, fix $R_1\leq R$, then for any $\rho\leq R_1$ there holds
$$
\omega(\nu_0 \rho)\leq\gamma(\rho)\omega(\rho)+R_{1}^{2}||f||_{\infty}
$$
and, therefore, for any integer $m>0$,
\begin{equation*}
\begin{split}
\omega(\nu_{0}^{m} R_1)&\leq\prod\limits_{k=0}^{m-1}\gamma(\nu_{0}^{k} R_1)\omega(R_1)+R_{1}^{2}||f||_{\infty}\left(\sum\limits_{k=1}^{m-1}\prod\limits_{i=1}^{k}\gamma(\nu_{0}^{m-i}R_1)+1\right)\\
&\leq\gamma(\nu_{0}^{m-1}R_1)^{m-1}\omega(R_1)+R_1^2||f||_{\infty}\frac{1-\gamma(\nu_{0}^{m-1}R_1)^{m-1}}{1-\gamma(\nu_{0}^{m-1}R_1)}
\end{split}
\end{equation*}
since $\gamma(r)$ is a decreasing function of $r$. Now, for any $r\leq R_1$ we can choose an integer $m>0$ such that $\nu_{0}^{m} R_1<r\leq \nu_{0}^{m-1} R_1$, and, therefore,
\begin{equation*}
\begin{split}
\omega(r)&\leq \gamma(r)^{m-1}\omega(R)+R_1^2||f||_{\infty}\frac{1-\gamma(r)^{m-1}}{1-\gamma(r)}\\
&\leq\frac{1}{\gamma(r)}\left(\frac{r}{R_1}\right)^{\ln(\gamma(r))/\ln(\nu_0)}+\frac{||f||_{\infty}}{1-\gamma(r)}R_1^2
\end{split}
\end{equation*}
Let $R_1=R^{1-\mu}r^{\mu}$ for $\mu\in(0,1)$ to obtain from above
$$
\omega(r)\leq\frac{1}{\gamma(r)}\left(\frac{r}{R}\right)^{(1-\mu)\ln(\gamma(r))/\ln(\nu_0)}+\frac{R^2||f||_{\infty}}{1-\gamma(r)}\left(\frac{r}{R}\right)^{2\mu}.
$$
From the definition of $\gamma(r)$ it is clear that $\gamma>1/2$, so the coefficient in front of the first term is bounded. Moreover, denote 
\begin{equation}\label{alpha_r}
\alpha(r)=\frac{(1-\mu)\ln(\gamma(r))}{\ln(\nu_0)}=\frac{(1-\mu)\ln(1-1/2C_{Har}(r))}{\ln(\nu_0)}
\end{equation}
and recall that $C_{Har}(r)\to\infty$ as $r\to 0$. It then follows that the estimate (\ref{pre_holder}) holds.
\end{proof}

We can now obtain continuity of weak solutions under a certain control of non doubling order of metric balls.
\begin{Th}\label{continuity_thm}
Let the estimate (\ref{pre_holder}) hold, and suppose that the non doubling order $\delta_y(r)$ satisfies (\ref{nondoub_order_control}) for all $y\in\Omega$, with $\lambda = (5\sigma-1)/(\sigma-1)$, $\sigma$ as in (\ref{sobolev}), and some $C>1$. Then every weak solution $u$ to (\ref{eq}) is continuous in $\Omega$. 
\end{Th}
\begin{proof}
First, note that if (\ref{nondoub_order_control}) holds with $\lambda = (5\sigma-1)/(\sigma-1)$, $\sigma$ as in (\ref{sobolev}), and $C$ as in definition of $C_{Har}$ (\ref{harnack}), then we have by (\ref{alpha_r})
\begin{equation*}
\begin{split}
\ln(r^{\alpha(r)})&=\frac{(1-\mu)\ln(1-1/2C_{Har}(r))}{\ln(\nu_0)}\ln r\\
&\leq -c\ln\left(1-\frac{1}{2C}\exp\left(-2\left[\frac{\nu_0r}{\delta(\nu_0 r)}\right]^{\frac{4\sigma}{\sigma-1}+1}\right)\right)\ln r
\to -\infty\quad\hbox{as}\quad r\to 0
\end{split}
\end{equation*}
and therefore, $r^{\alpha(r)}\to 0$ as $r\to 0$ . Next, it is easy to see, that for any $C_1,C_2\geq 1$ and $f(r)<1$, $f(r)\to 0$ as $r\to 0$ there holds
$$
\lim_{r\to 0}\frac{\ln(1-\frac{f(r)}{C_1})}{\ln(1-\frac{f(r)}{C_2})}=const.
$$
Therefore, if (\ref{nondoub_order_control}) holds with some $C\geq 1$, it holds with any $C\geq 1$ and we get that the first term on the right in (\ref{pre_holder}) converges to zero as $r\to 0$.\\
To estimate the second term use $2C_{Har}r^{2\mu}=r^{2\mu}/(1-\gamma(r))$ and write
$$
\ln\left(\frac{r^{2\mu}}{1-\gamma(r)}\right)=\frac{1}{\ln(\gamma(r))}\left(-\ln(\gamma(r))\ln(1-\gamma(r))+2\mu\ln(r)\ln(\gamma(r))\right).
$$
Since $\gamma(r)\to 1$ as $r\to 0$, the first term in parentheses on the right hand side converges to $0$. By (\ref{nondoub_order_control}) and the definition of $\gamma(r)$ the second term approaches positive infinity. Since $\lim_{r\to 0^{+}}\ln(\gamma(r))=0^{-}$ we therefore obtain
$$
\lim_{r\to 0^{+}}\ln\left(\frac{r^{2\mu}}{1-\gamma(r)}\right)=-\infty
$$ 
or, equivalently,
$$
\lim_{r\to 0^{+}}\left(2C_{Har}(r)r^{2\mu}\right)=\lim_{r\to 0^{+}}\left(\frac{r^{2\mu}}{1-\gamma(r)}\right)=0.
$$
It then follows that for any $x,x'\in B(y, R)$ there holds
$$
\lim_{x\to x'}|u(x)-u(x')|=0
$$
\end{proof}

\section{Subunit metric}\label{Subunit metric}
In this section we describe a certain metric associated to the operator, a subunit metric. This is the metric used by Fefferman and Phong to give their criterion of subellipticity \cite{Fef}. Although it is not required in Theorem \ref{main_intro} that the metric is subunit, it will be clear that this is one natural choice of the metric, for which the requirements of the theorem may be satisfied. In particular, we show that some of the assumptions of Theorem \ref{main_intro}, such as the existence of the accumulating sequence of cutoff functions, do indeed hold true for subunit metrics associated to a certain class of degenerate operators.

\subsection{Definitions}
Let $Q(x,\xi)$ be nonnegative semidefinite continuous quadratic form, $\Omega\subset\RR^n$ an open bounded domain.

\begin{Def}
A Lipschitz curve $\gamma\;:\;[0,r]\rightarrow\Omega$ is subunit with respect to $Q(x,\xi)$ if
\begin{equation*}
(\gamma'(t)\cdot\xi)^2\leq Q(\gamma(t),\xi)\  a.e.\;t\in[0,r],\;\xi\in\mathbb{R}^{n}.
\end{equation*}
\end{Def}

\begin{Def}
We define
\begin{equation}\label{subunit_dist}
d(x,y)=\inf\{r>0\;:\;\gamma(0)=x,\;\gamma(r)=y,\;\gamma\hbox{ is subunit in $\Omega$}\}.
\end{equation}
\end{Def}

\begin{Def}
If $d$ is finite on $\Omega\times\Omega$ we define subunit balls $B(x,r)$ by
\begin{equation*}
B(x,r)=\{y\in\Omega\;:\;d(x,y)<r\},\  \  x\in\Omega,\;0<r<\infty.
\end{equation*}
\end{Def}
We will denote $E(x,r)$ the usual Euclidean ball
\begin{equation*}
E(x,r)=\{y\in\Omega\;:\;|x-y|<r\},\  \  x\in\Omega,\;0<r<\infty.
\end{equation*}

First, note that the boundedness condition (\ref{upper}) is almost immediate. That is, there exists a constant $C>0$ s.t. for any $x\in\Omega$, any $r>0$ s.t. $E(x,Cr)\subset\Omega$ there holds
$B(x,r)\subset E(x,Cr)$. Indeed, note that we can define the Euclidean distance $|x-y|$ between two points $x,y\in\Omega$ as follows
\begin{equation*}
|x-y|=\inf\{r>0\;:\;\gamma(0)=x,\;\gamma(r)=y,\;\gamma\hbox{ is a rectifiable curve in $\Omega$}\}.
\end{equation*}
It is clear from (\ref{subunit_dist}) that $|x-y|\leq d(x,y)$ since the set of all rectifiable curves is wider than the set of all subunit curves, and thus (\ref{upper}) follows with $C=1$.

The containment condition (\ref{contain}) will be assumed to hold for subunit balls, namely, $\forall x\in\Omega\;\exists\alpha_x:\RR_{+}\rightarrow\RR_{+}$ increasing and such that $\alpha_x(r)>0$ for any $r>0$ and
\begin{equation*}
E(x,\alpha_x(r))\subset B(x,r),\  \forall x\in\Omega,\  \forall r>0\;\hbox{s.t.}\;B(x,r)\subset\Omega
\end{equation*}


\subsection{Existence of Cutoff Functions}\label{section_cutoff_functions}
The following lemma shows the existence of a special cutoff function required in the proof of the weak type logarithmic estimates, Lemma \ref{lemma_weak1}, in the case of a subunit metric defined as above.
\begin{Lemma}\label{spec_cutoff_lemma}
Let $Q(x,\xi)=\xi'Q(x)\xi$ be a continuous nonnegative semidefinite quadratic form. Suppose that the subunit metrics associated to $Q(x)$ satisfies the containment condition (\ref{contain}) and the boundedness condition (\ref{upper}). Then for each ball $B(y,r)$ with $y\in\Omega$, $0<r<\eta\, dist(y,\partial\Omega)$ there exists a cutoff function $\phi_{r}\in Lip(\Omega)$ satisfying (\ref{spec_cutoff}), namely
\begin{equation*}
\begin{cases}
\supp(\phi_{r}) &\subseteq B(y,r+\delta(r))\\
\{x:\phi_{r}(x)=1\} &\supseteq B(y,r+\delta(r)/2)\\
||[\nabla\phi_{r}]_Q||_{L^{\infty}(B(y,r)} &\leq \frac{C}{\delta(r)}
\end{cases}
\end{equation*}
where $\delta(r)=\delta_y(r)$ as in (\ref{nondoub_order}).
\end{Lemma}

\begin{proof}
For any $\varepsilon\geq 0$ let $Q^{\varepsilon}(x,\xi)=Q(x,\xi)+\varepsilon^2|\xi|^2$. Then as it has been shown in \cite[Lemma 66]{Saw} the subunit metric $d^{\varepsilon}(x,y)$ associated to $Q^{\varepsilon}$ satisfies
\begin{equation*}
[\nabla_{x}d^{\varepsilon}(x,y)]_{Q}\leq \sqrt{n}, \quad x,y\in\Omega
\end{equation*}
uniformly in $\varepsilon>0$. Moreover, $d^\varepsilon(\cdot,y)\nearrow d(\cdot,y)$, the convergence is monotone and $d$ is continuous, therefore, $d^\varepsilon(\cdot,y)\to d(\cdot,y)$ uniformly on compact subsets of $\Omega$. Define $g(t)$ to vanish for $t\geq r_2$, to equal $1$ for $t\leq r_1$ and to be linear on the interval $[r_1,r_2]$, where $r_1=r$, $r_2=r+\delta(r)/2$. Let $\phi_r(x)=g(d^{\varepsilon^*}(x,y))$, with $\varepsilon^*$ to be chosen later. We have
$\phi_r(x)=1$ when $d^{\varepsilon^*}(x,y)\leq r$, and since $d^{\varepsilon^*}(x,y)\leq d(x,y)$ we have $\phi_r(x)=1$ when $d(x,y)\leq r+\delta(r)/2$. Next, $\phi_r(x)=0$ when $d^{\varepsilon^*}(x,y)>r+\delta(r)$ and by choosing $\varepsilon^{*}$ small enough, we obtain that $\phi_r(x)=0$ when $d(x,y)>r+\delta(r)$. This shows that $\supp(\phi_{r}) \subseteq B(y,r+\delta(r))$ and $\{x:\phi_{r}(x)=1\} \supseteq B(y,r+\delta(r)/2)$.
Next,
\begin{equation}\label{cut_1}
[\nabla\phi_r(x)]_{Q}\leq ||g'||_{\infty}[\nabla d^{\varepsilon^{*}}]_{Q}\leq\frac{C}{\delta(r)}\sqrt{n}.
\end{equation}
This completes the proof.
\end{proof}

The following lemma is a variant of Proposition 68 in \cite{Saw}.
\begin{Lemma}\label{cutoff_lemma}
Let $Q(x,\xi)=\xi'Q(x)\xi$ be continuous nonnegative semidefinite quadratic form. Suppose that the subunit metrics associated to $Q(x)$ satisfies (\ref{contain}), (\ref{upper}).
Then there exists an accumulating sequence of cutoff functions satisfying (\ref{cutoff}), namely, for any $\nu<1$,
\begin{equation*}
\begin{cases}
E_1&=\supp(\psi_1) \subset B(y,r)\\
B(y,\nu r) &\subset \{x:\psi_j(x)=1\},\;j\geq 1\\
E_{j+1}=\supp(\psi_{j+1}) &\Subset \{x:\psi_j(x)=1\},\;j\geq 1\\
\frac{\displaystyle |E_j|}{\displaystyle|E_{j+1}|}&\leq C,\;j\geq 1\\
\psi_j &\hbox{is Lipschitz},\;j\geq 1\\
||[\nabla\psi_j]_Q||_{L^{\infty}(B(y,r))} &\leq \frac{C}{\displaystyle (1-\nu)\delta(\nu r)\left(1-\frac{\delta(\nu r)}{r}\right)^{j}},\;j\geq 1
\end{cases}
\end{equation*}
where $\delta(r)=\delta_y(r)$ as in (\ref{nondoub_order}).
\end{Lemma}

\begin{proof}
The sequence will be constructed in a similar way as a function $\phi_r(x)$ in Lemma \ref{spec_cutoff_lemma}. As in the proof of Lemma \ref{spec_cutoff_lemma} for any $\varepsilon>0$ let $Q^{\varepsilon}(x,\xi)=Q(x,\xi)+\varepsilon^2|\xi|^2$. Then
\begin{equation*}
[\nabla_{x}d^{\varepsilon}(x,y)]_{Q}\leq \sqrt{n}, \quad x,y\in\Omega
\end{equation*}
and $d^\varepsilon(\cdot,y)\to d(\cdot,y)$ uniformly on compact subsets of $\Omega$.
Define $\phi_j(t)$, $j\geq 1$, to vanish for $t\geq r_j$, to equal $1$ for $t\leq r_{j+1}$ and to be linear on the interval $[r_{j+1},r_{j}]$. Here the radii $r_j$ are defined as follows
\begin{equation*}
\begin{split}
r_1&=r-(1-\nu)\delta(\nu r)\\
r_j&=r-(1-\nu)\delta(\nu r)\sum\limits_{i=0}^{j-1}\left(1-\frac{\delta(\nu r)}{r}\right)^i,\;\;\forall j\geq 2
\end{split}
\end{equation*}
Now, let $\psi_j(x)=\phi_{j}(d^{\varepsilon_{j}}(x,y))$, with $\varepsilon_j$, decreasing in $j$, to be chosen later. We have
$\psi_{j}=1$ on $\{x:d^{\varepsilon_{j}}(x,y)\leq r_{j+1}\}$, and $\psi_j=0$ on $\{x:d^{\varepsilon_{j}}(x,y)\geq r_{j}\}$. Moreover, $d\geq d^{\varepsilon_{j+1}}\geq d^{\varepsilon_{j}}$ and therefore, $\psi_{j}=1$ on $E_{j+1}:=\supp(\psi_{j+1})=\{x:d^{\varepsilon_{j+1}}(x,y)\leq r_{j+1}\}$, and $E_{j+1}\supseteq B(y,r_{j+1})$. 
Choosing $\varepsilon_j$ small enough we also obtain that $\psi_j=0$ on $\{x:d(x,y)\geq r_j+\nu(r_j-r_{j+1})\}=B(y,r_j+\nu(r_j-r_{j+1}))^{c}$, and thus $E_j=\supp(\psi_j)\subseteq B(y,r_j+\nu(r_j-r_{j+1}))$ with $\psi_{j}=1$ on $\supp(\psi_{j+1})$.
Next, denoting $r_{\infty}=\lim_{j\rightarrow\infty}r_j$ and writing
$$
r-r_{\infty}=\sum\limits_{j=0}^{\infty}(r_j-r_{j+1})=\sum\limits_{j=0}^{\infty}(1-\nu)\delta(\nu r)\left(1-\frac{\delta(\nu r)}{r}\right)^{j}=r(1-\nu)
$$
one can see that $r_{\infty}=\nu r$ and therefore 
\begin{equation*}
B(y,\nu r)\subset E_j,\;\;\forall j\geq 1.
\end{equation*}
Moreover
\begin{equation}\label{smallness}
r_j-r_{j+1}=(1-\nu)\delta(\nu r)\left(1-\frac{\delta(\nu r)}{r}\right)^{j}\leq(1-\nu) \delta(\nu r)\leq(1-\nu) \delta(r_j)\leq\delta(r_j)
\end{equation}
We can now easily see that the fourth condition in (\ref{cutoff}) is satisfied. Indeed, using \ref{smallness} we have
\begin{equation*}
\begin{split}
\frac{|E_{j}|}{|E_{j+1}|}&\leq\frac{|B(y,r_j+\nu(r_j-r_{j+1}))|}{|B(y,r_{j+1}))|}
=\frac{|B(y,r_{j+1}+(1+\nu)(r_{j}-r_{j+1}))|}{|B(y,r_{j+1})|}\\
&\leq \frac{|B(y,r_{j+1}+\delta(r_{j+1}))|}{|B(y,r_{j+1})|}\leq C
\end{split}
\end{equation*} 
by the definition of $\delta(r)$. To verify the last condition in (\ref{cutoff}) write
$$
[\nabla\psi_j]_{Q}\leq ||\nabla\phi_{j}||_{\infty}[\nabla d^{\varepsilon_{j}}]_{Q}\leq\frac{C}{r_{j}-r_{j+1}}\sqrt{n}
\leq\frac{C}{(1-\nu)\delta(\nu r)\left(1-\frac{\delta(\nu r)}{r}\right)^{j}}
$$
Note also, that since $\delta(r)/r\to 0$ as $r\to 0$, the expression $(1-\delta(\nu r)/r)$ is bounded from below by $(1-\delta(\nu R_0)/R_0)>0$ for all $r\leq R_0$.
This concludes the proof of the existence of the sequence satisfying (\ref{cutoff}).
\end{proof}

\begin{Rem}
Note that the condition that $Q(x)$ is continuous cannot be easily omitted. In \cite{Zhong} the author constructs an example of a discontinuous solution to a degenerate linear elliptic equation (see Theorem 1.3 and Conjecture 6). However, the matrix $Q$ in that case is discontinuous and this requirement seems to be essential for the construction.
\end{Rem}


\section{Model examples}\label{examples}
In this section we construct some model examples of operators such that the subunit metrics associated to them define nondoubling balls satisfying (\ref{nondoub_order_control}). We start with a linear example, which can then be easily modified into a nonlinear one.
\subsection{Linear Example}
Consider a model example of a linear infinitely degenerate second order operator in two dimensions $L=\nabla^{T}A(x,y)\nabla$, where the matrix $A$ is defined as follows
\begin{equation}\label{model}
A(x,y)=\left(
\begin{array}{cc}
1& 0\\
0& f(x)^2
\end{array}
\right)
\end{equation}
Here the function $f(x)\in C^{\infty}(\RR)$ is even, $f(x)>0,\,\forall x\neq 0$, and $f(x)$ can vanish at $x=0$. We will also assume that $f$ is an increasing function on $\RR_{+}$.
We would like to consider a function $f$ which vanishes at $x=0$ together with all its derivatives. This will imply that the operator $L$ degenerates to infinite order on the $y$-axis.
First, we show that the subunit metric balls defined by this operator satisfy the containment condition (\ref{contain}). More precisely, we give the estimate on the size of subunit balls in terms of boxes. The following lemma deals with subunit balls that have centers on the $y$-axis. These balls are the ``smallest'' among the balls of the fixed radius.

\begin{Lemma}\label{boxes}
Let the operator $L$ be defined by the matrix $A(x,y)$ as in (\ref{model}). Define the box $Q_r(x,y)=[x-r,x+r]\times [y-rf(r/2),y+rf(r/2)]$, and $\tilde{Q}_r(x,y)=[x+r/2,x+3r/4]\times[y-r/4f(r/2),y+r/4f(r/2)]$. Then for any $r>0$ the subunit balls $B_r=B((0,y),r)$ satisfy
\begin{equation}\label{boxes1}
\tilde{Q}_r(0,y)\subset B_r \subset Q_{r}(0,y)
\end{equation}
and consequently, 
\begin{equation}\label{boxes_size}
r^2/8f(r/2)\leq |B_r|\leq 4r^2f(r/2)
\end{equation}
\end{Lemma}
\begin{proof}
Since the operator matrix only depends on the first variable, it is enough to prove the statement for $y=0$. To prove the first inclusion in (\ref{boxes1}), first note, that any horizontal line segment is a subunit curve. Next, consider a point $p=(x_0,r/4f(r/2))$ on the ``top side'' of $\tilde{Q}(0,0)$, where $r/2\leq x_0\leq 3r/4$. Let a subunit curve $\gamma=\gamma_1\cup\gamma_2$ connect the origin to the point $p$. Here, $\gamma_1$ is a horizontal segment, connecting $(0,0)$ to $(x_0,0)$ and $\gamma_2(t)$ is defined as follows
\begin{equation*}
\gamma_{2}(t)=\left(x_0,f(r/2)t\right),\;\;t\in\left[0,\frac{r}{4}\right].
\end{equation*}
Since $f$ is an increasing function on $\RR_{+}$, it is easy to check that $(\gamma'(t)\cdot\xi)^2\leq \xi^TA(\gamma(t))\xi$ for any $\xi\in\RR^n$, and therefore, 
$$
d((0,0),p)\leq |x_0-\frac{r}{2}|+\frac{r}{4}<r.
$$ 
Therefore, $p\in B(0,r)$. Moreover, it is clear that any other point in $\tilde{Q}_{r}$ can be connected to the origin by a similarly constructed curve, so that the distance to the origin is less than $r$. This concludes the proof that $\tilde{Q}_{r}\in B_r$.
To show  the other inclusion, let $\gamma(t)$ be the minimizing curve connecting the origin to any point on the boundary $\partial Q_{r}$. First, let the point $(x,y)$ belong to the top or the bottom edge of $\partial Q_{r}$, i.e. $|y|=rf(r/2)$. Without loss of generality we can also assume, $x\geq 0$. The curve $\gamma(t)$ is thus a subunit curve satisfying
$$
\gamma(0)=(0,0),\;\;\gamma(T)=(x,y),\;\;T=d((0,0),(x,y)).
$$
Then we have
\begin{equation}\label{T}
rf(r/2)=|y-0|=\left|\int\limits_{0}^{T}\gamma_2'(t)dt\right|\leq\int\limits_{0}^{T}\left|\gamma_2'(t)\right|dt
\leq\int\limits_{0}^{T}f(\gamma_1(t))dt
\end{equation}
To estimate $|\gamma_1(t)|$ first note the following
$$
T=d((0,0),(0,y))=d((0,0),(\gamma_1(t),\gamma_2(t)))+d((\gamma_1(t),\gamma_2(t)),(0,y))
$$
Moreover, since the Euclidean distance can be defined in the same way as (\ref{subunit_dist}) but without the restriction on the curve being subunit, we have $|x-y|\leq d(x,y)$ $\forall x,y\in\Omega$. Thus, we obtain
$$
T=d((0,0),(0,y))\geq \sqrt{\gamma_1(t)^2+\gamma_2(t)^2}+\sqrt{\gamma_1(t)^2+(y-\gamma_2(t))^2}\geq 2|\gamma_1(t)|
$$
or $|\gamma_1(t)|\leq T/2$ and therefore from (\ref{T})
$rf(r/2)\leq Tf(T/2)$.
By assumption, the function $xf(x)$ is strictly increasing for $x>0$ and thus $T\geq r$. Now, if the point $(x,y)\in \partial Q_{r}$ satisfies $|x|=r$ it is obvious that $d((0,0),(x,y))\geq r$. This completes the proof.
\end{proof}

We now would like to define a function $f(x)$ such that the subunit balls associated to the operator (\ref{model}) are non doubling, but the non doubling order satisfies the estimate (\ref{nondoub_order_control}). The model (\ref{model}) will then fall under the framework of Theorem \ref{main_intro}. The next proposition provides an example of such function $f$.

\begin{Prop}
Let the function $h(x)$ be defined as follows
\begin{equation*}
h(x)=\left(-\frac{1}{\ln\left(1-\exp\left(\frac{1}{(\ln x)^{1/3}}\right)\right)}\right)^{\frac{1}{\lambda}}
\end{equation*}
with $\lambda = (5\sigma-1)/(\sigma-1)$, $\sigma$ as in (\ref{sobolev}). Define 
$$
f(x)=\exp\left(-\int\limits_{x}^{1}\frac{1}{th(t)}dt\right)
$$ 
and let the operator $L$ be given by the matrix (\ref{model}). Then the subunit balls associated to $A(x,y)$ with centers on the $y$-axis are non doubling, and the non doubling order satisfies (\ref{nondoub_order_control}).
\end{Prop}
\begin{proof}
First, it is easy to check that $\lim_{x\to 0^{+}}h(x)=0^{+}$. Moreover, $\lim_{x\to 0^+}\int\limits_{x}^{1}1/th(t)dt=+\infty$, and therefore $f(0)=0$. To see that $f(x)$ is infinitely degenerate at the origin we first show that the subunit balls associated to $A(x,y)$ with centers on the $y$-axis are non doubling.
To show that the subunit balls are non doubling according Lemma \ref{boxes} it is enough to consider the ratio $f(r)/f(r/2)$. Indeed, using (\ref{boxes_size}) we have
\begin{equation}\label{model_nondoub}
\begin{split}
\frac{|B((0,y),2r)|}{|B((0,y),r)|}&\geq c \frac{f(r)}{f(r/2)}=c \exp\left(\int\limits_{r/2}^{r}\frac{1}{th(t)}dt\right)\\
&\geq c\exp\left(\frac{1}{rh(r)}\int\limits_{r/2}^{r}dt\right)=c\exp\left(\frac{1}{2h(r)}\right)\to\infty,\hbox{\quad as \quad}r\to 0
\end{split}
\end{equation}
On the other hand, let $\delta(r)=rh(r/2)$, then
\begin{equation*}
\begin{split}
\frac{|B((0,y),r+\delta(r))|}{|B((0,y),r)|}&\leq C \frac{f((r+\delta(r))/2)}{f(r/2)}=C \exp\left(\int\limits_{r/2}^{r/2+rh(r/2)/2}\frac{1}{th(t)}dt\right)\\
&\leq C\exp\left(\frac{1}{rh(r/2)/2}\int\limits_{r/2}^{r/2+rh(r/2)/2}dt\right)=Ce
\end{split}
\end{equation*}
which implies that the non doubling order $\delta(r)\approx rh(r/2)$. It is a simple exercise to verify that $\delta(r)=rh(r/2)$ satisfies condition (\ref{nondoub_order_control}). Finally, to see that the function $f(x)$ has infinite degeneracy at $x=0$, let us assume the opposite, there exists an integer $N$ such that $f^{(N)}(0)\neq 0$. Then it follows that $f(r)\approx r^N$ as $r\to 0$ and
$$
\frac{f(2r)}{f(r)}\leq C,\hbox{\quad as\quad} r\to 0
$$
which contradicts (\ref{model_nondoub}).

\end{proof}
\subsection{Nonlinear Example}

First, note that although the matrix $A(x,y)$ from the previous example is linear and $C^{\infty}$, the proof above does not rely on differentiability properties of $A(x,y)$. Now, let the matrix $A$ depend also on the solution, $A(x,y,u)$. Recall that all the assumptions of the regularity result, Theorem \ref{main_intro}, are stated for the matrix $Q$ from the structural assumption (\ref{struc})
$$
k\, \xi^TQ(x,y)\xi \leq \xi^TA(x,y,z)\xi\leq K\,\xi^TQ(x,y)\xi
$$
for a.e. $x,y\in\Omega$ and all $z\in\mathbb{R}$, $\xi\in\mathbb{R}^n$. We therefore can modify the previous example as follows. Let 
\begin{equation}\label{model_nonlin}
A(x,y,u)=\left(
\begin{array}{cc}
1& 0\\
0& \phi(u)f(x)^2
\end{array}
\right)
\end{equation}
where $\phi(t):\RR\to\RR$ is a bounded measurable function and $c\leq \phi(t)\leq C$ for all $t\in\RR$ and some positive constants $c$ and $C$. Then obviously, for all $x,y$ and all $z\in\mathbb{R}$, $\xi\in\mathbb{R}^n$ there holds
$$
c\, \xi^TQ(x,y)\xi \leq \xi^TA(x,y,z)\xi\leq C\,\xi^TQ(x,y)\xi
$$
with
\begin{equation*}
Q(x,y)=\left(
\begin{array}{cc}
1& 0\\
0& f(x)^2
\end{array}
\right)
\end{equation*}
and therefore, the previous example applies.

\bibliographystyle{amsplain}

\bibliography{library2}

\end{document}